\documentstyle[12pt]{article}
\begin{document}
\newtheorem{lem}{Lemma}[section]
\newtheorem{rem}{Remark}[section]
\newtheorem{question}{Question}[section]
\newtheorem{prop}{Proposition}[section]
\newtheorem{cor}{Corollary}[section]
\newtheorem{thm}{Theorem}[section]
\newtheorem{definition}{Definition}[section]
\newtheorem{exam}{Example}[section]
\newtheorem{conj}{Conjecture}[section]
\newtheorem{openproblem}{Open Problem}[section]
\def\av{{\int \hspace{-2.25ex}-} }
\title
{Some remarks on singular solutions of
 nonlinear elliptic
equations. I. }
\author{Luis Caffarelli\thanks{Partially
 supported by
        NSF grant DMS-0654267}\\
Department of Mathematics
\\
The University of Texas
\\
Austin, TX 78712
\\
\\
YanYan Li\thanks{Partially
 supported by
        NSF grant DMS-0401118 and DMS-0701545.}\\
Department of Mathematics\\
Rutgers University\\
110 Frelinghuysen Road\\
Piscataway, NJ 08854\\
\\
Louis Nirenberg\thanks{Not even partial support.}
 \\
Courant Institute\\
251 Mercer Street \\
New York, NY 10012\\
}

\date{ }
\maketitle

\centerline{  {\it Dedicated to Felix Browder on his $80^{\mbox{th}}$
 birthday}}

\medskip

%%%%%  The following two lines should be changed
%\input amstex
\input { amssym.def}
%%%%%%%%%

\setcounter {section} {0}

\section{Introduction}

\subsection{}
Understanding the behavior of solutions of
partial differential equations near an isolated
singularity is of basic importance in the study
of partial differential equations.
A classical theorem of
B\^ocher (see e.g.
\cite{Bo}) states that a positive harmonic function $u$ in a
punctured ball $B_1\setminus\{0\}$ in $\Bbb R^n$  must be of the form
$$
u(x)=
\left\{
\begin{array}{rl}
-a\log |x|+h(x), & \mbox{if}\ n=2,\\
a|x|^{2-n}+h(x), & \mbox{if}\ n\ge 3,
\end{array}
\right.
$$
where $a$ is a nonnegative constant and $h$ is
a harmonic function in the ball $B_1$.

For nonlinear partial differential equations or systems,
it might be difficult to give as complete a description
of the behavior of solutions near an isolated singularity.
On the other hand,  some studies  only require
certain partial information on  the behavior of solutions near an isolated
singularity.  Part of this paper, e.g. Theorem
\ref{thm1} and Theorem \ref{thmnew1}, address
some  issues which arise naturally in
applications of the method
of moving planes.

Let $F\in C^1(\Omega\times \Bbb R\times \Bbb R^n\times {\cal S}^{n\times n})$,
where ${\cal S}^{n\times n}$ denotes the set of $n\times n$
real symmetric  matrices and $\Omega$ is a
domain (bounded connected
open set)  in the $n-$dimensional
Euclidean space $\Bbb R^n$.
Without loss of generality we assume that $\Omega$ contains the origin $0$.
Throughout the paper we use $B_r(x)$ to denote
a ball  of radius $r$ and centered at $x$,
and use $B_r$ to denote $B_r(0)$.

We assume that $F(x, v, \nabla v, \nabla^2 v)
$ is an elliptic, but not necessarily uniformly 
elliptic,  operator:
\begin{equation}
\left( -\frac{\partial F}{\partial M_{ij} }(x, s, p, M)\right)
>0\qquad \forall\ (x, s, p, M)\in
\Omega \times \Bbb R\times \Bbb R^n\times {\cal S}^{n\times n}.
\label{elliptic}
\end{equation}

In our first theorem, the function $u\in
C^2(\Omega\setminus \{0\})$ has the following property
for some $\bar r>0$:
\begin{equation}
\mbox{For any}\ V\in \Bbb R^n, w(x):= u(x)+V\cdot x\
\mbox{satisfies}\
\inf_{B_r\setminus \{0\}}w=\min_{ \partial B_r }w\
\forall\ 0<r<\bar r.
\label{superaffine}
\end{equation}
In dimension $n\ge 2$,
a superharmonic function $u\in C^0(B_3\setminus \{0\})$
 satisfying
$\inf_{B_3\setminus \{0\}}u>-\infty$ has  the above
property. On the other hand,
the function \begin{equation}
u(x):=
\left\{
\begin{array}{rl}
4, &  -3<x<0,\\
1,&0\le x<3,
\end{array}
\right.
\label{e1}
\end{equation}
a harmonic function in $(-3, 3)\setminus \{0\}$,
does not satisfy condition (\ref{superaffine})
for any $\bar r>0$.
To see this, simply take $w(x)=u(x)+ x$.
We start with a result related to the strong maximum principle.
\begin{thm}  For $n\ge 1$, let $\Omega$ be a domain
 in $\Bbb R^n$ containing the origin
$\{0\}$, and  let $F\in C^1(\Omega\times \Bbb R\times
 \Bbb R^n\times {\cal S}^{n\times n})$
  satisfy
(\ref{elliptic}).   Assume that   $u\in
C^2(\Omega\setminus \{0\})$ satisfies (\ref{superaffine})
for some $\bar r>0$,
$v\in C^2(\Omega)$,
$$
u>v\qquad \mbox{in}\ \Omega\setminus \{0\},
$$
$$
F(x, u, \nabla u, \nabla^2 u)\ge F(x, v, \nabla v, \nabla^2v) \qquad
\mbox{in}\ \Omega\setminus \{0\}.
$$
 Then
\begin{equation}
\liminf_{ |x|\to 0}(u-v)(x)>0.
\label{aa1}
\end{equation}
\label{thm1}
\end{thm}

\begin{rem} In dimension $n\ge 2$,  condition  (\ref{superaffine})
in the above theorem
can be replaced by 
 $a_{ij}(x)\partial_{ij}u\le
C$ 
in $B_{\bar r}\setminus \{0\}$
for some $\bar r>0$, where $(a_{ij}(x))$
  is some 
 positive definite matrix functions which is
continuous in $B_{\bar r}$ if $n\ge 3$ and
H\"older continuous in $B_{\bar r}$ if $n=2$, 
and $C$ is some constant.  Indeed, if $C=0$, then $u$ satisfies
(\ref{superaffine}); if $C\ne 0$, we may reduce it to
the $C=0$ case by working with $u(x)- A |x|^2$ and $
v(x)-A |x|^2$, for some large constant $A$, instead of $u$ and $v$,
and working with the new $F$. 
On the other hand, it is not the
case for $n=1$.  This can be seen by taking
$F(x,u,\nabla u, \nabla^2 u)=-u''$, $v(x)=1-x$ and $u$
as  in (\ref{e1}).
\end{rem}

\begin{rem}
If the assumption ``$u>v$ in
$\Omega\setminus \{0\}$'' in  the above theorem
is replaced by ``$u \ge v$ in
$\Omega\setminus \{0\}$'', then, in view
of the strong maximum principle,
either $u\equiv v$ in $\Omega\setminus\{0\}$ or
(\ref{aa1}) holds.
\label{rembb1}
\end{rem}

Theorem \ref{thm1}  is related to the following: Consider a viscosity
 supersolution $u$ 
 of a fully
nonlinear elliptic equation outside a singular point $\{0\}$
(or a singular set $S$).  Under what condition
is $u$  also a supersolution
across $\{0\}$ (or across $S$)?  We plan to address this issue in
subsequent papers.

 Condition
 (\ref{superaffine}) is
related to the following
\begin{definition}
A function $u\in C^0(\Omega)$ is called superaffine in a domain
$\Omega$ if in any subdomain $D$
 ($u+$ any linear function) achieves
its minimum in $\overline D$ on $\partial D$.
\label{def1}
\end{definition}
This notion was also introduced by Harvey
and Lawson in  \cite{HL}  at the same time.
For a $u\in C^2(\Omega)$,  $u$ is
 superaffine in $\Omega$  if and only if
 the lowest eigenvalue $\lambda_1(\nabla^2 u)$
of the Hessian of $u$ is $\le 0$ in $\Omega$ --- as
is easily seen (see the proof of Proposition \ref{prop1} below).

Throughout this paper we always
order  the eigenvalues of the Hessian of
 a $C^2$ function $u$ by
$$
\lambda_1(\nabla^2 u)\le \lambda_2(\nabla^2 u)\le
\cdots\le \lambda_n(\nabla^2 u).
$$

Existence and uniqueness of continuous solutions to
general degenerate equations,  including
$$
\lambda_p(\nabla^2 u)+\cdots+\lambda_{p+q}(\nabla^2 u)=0
$$
for fixed $p$ and $q$ as examples, have been
obtained in  \cite{HL}.

Condition  (\ref{superaffine})  is slightly stronger than
$u$ being superaffine in $\Omega\setminus\{0\}$, for example
it is not satisfied by $|x|^a$, $0<a<1$, which is
superaffine in $\Omega\setminus \{0\}$.
However the following holds.

\begin{lem} For $n\ge 2$,
let $u\in C^2(B_1\setminus \{0\})$ satisfy
\begin{equation}
\lambda_1(\nabla^2 u)+  \lambda_2(\nabla^2 u)\le
0,\qquad \mbox{in}\ B_1\setminus \{0\},
\label{mn1}
\end{equation}
and
\begin{equation}
\inf_{ B_1\setminus \{0\} }u>-\infty.
\label{abab1}
\end{equation}
Then
condition  (\ref{superaffine}), with $\bar r=1$, holds.
\label{lem0-1}
\end{lem}

\begin{rem} In Lemma \ref{lem0-1}, condition (\ref{mn1})
cannot be replaced by the following    weaker
condition: For some $0<a<1$,
\begin{equation}
 \lambda_1(\nabla^2 u)+ a \lambda_2(\nabla^2 u)\le 0,\qquad
\qquad \mbox{in}\ B_1\setminus \{0\}.
\label{mn2}
\end{equation}
Indeed,
 the function
$u(x)=|x|^{1-a}$
 satisfies (\ref{mn2}) (with the equality).
\end{rem}

Note that
$$
\lambda_1(\nabla^2 u)+ a \lambda_2(\nabla^2 u)\le
\frac {1+a}2 (\lambda_1(\nabla^2 u)+  \lambda_2(\nabla^2 u)),
$$
for  $0\le a<1$.
So (\ref{mn2}) is weaker than (\ref{mn1}).

\begin{rem}
The assumption (\ref{abab1}) in  Lemma \ref{lem0-1}
 cannot be removed, as can be seen by
taking
$$
u(x)=
\left\{
\begin{array}{rl}
\log|x|,& n=2,\\
-|x|^{2-n},& n\ge 3.
\end{array}
\right.
$$
\end{rem}

Theorem \ref{thm1}  would not hold if condition (\ref{superaffine})
 were dropped.  Here are some examples.

\medskip

\noindent{\bf Example 1.}\ For$n\ge 1$,  $l=2, 3, 4, ...$, let
$\alpha=\frac{2l-2}{2l-1}$.  Clearly $\alpha \in (0, 1)$ and
$2l=\frac {2-\alpha}{1-\alpha}$.  For $u(x):=|x|^\alpha$ and
$v\equiv 0$, we have $u>v$ in $B_1\setminus\{0\}$ and
$$
F(\nabla u, \nabla^2 u):= -\Delta u+ (n+\alpha-2))\alpha^{\frac
1{\alpha-1}}|\nabla u|^{2l}=0=F(\nabla v, \nabla^2 v)\qquad
\mbox{in}\ B_1\setminus\{0\}.
$$
But
\begin{equation}
\lim_{ |x|\to 0}(u-v)(x)=0.
\label{h}
\end{equation}

\medskip

\noindent{\bf Example 2.}\ For $n=2$, $0<\alpha<1$,
$\epsilon= 1-\sqrt{\alpha} \in (0, 1)$,  let
$$
b_1(x)=-(1-\epsilon)\frac {x_2}{|x|}, \ \
b_2(x)=(1-\epsilon)\frac{ x_1}{|x|}, \ \
a_{ij}(x)=\delta_{ij}-b_i(x)b_j(x).
$$
Consider $u(x):=|x|^\alpha$, $v\equiv 0$, we have
 $u>v$ in $B_1\setminus\{0\}$ and
\begin{eqnarray*}
F(x, \nabla^2 u)&:=&- a_{ij}(x)u_{ij}(x)
=-\alpha |x|^{\alpha-2}
a_{ij}(x)
\left( \delta_{ij}-(2-\alpha)\frac {x_ix_j}{  |x|^2 }\right)
\\
&=&
-(\alpha-(1-\epsilon)^2)\alpha |x|^{\alpha-2}
= 0=F(x, \nabla^2 v)\qquad\qquad \mbox{in}\ B_1\setminus\{0\}.
\end{eqnarray*}
But (\ref{h}) holds.

\medskip

\noindent{\bf Example 3.}\ For $n\ge 2$, $0<\alpha<1$,
$ \epsilon=
(1-\alpha)/(n-1)  \in (0, 1)$,
  let $O(x)$ be an $n\times n$
orthogonal matrix functions in $L^\infty(B_1\setminus \{0\})$ such that
$$
O(x)_{ij}\frac {x_j}{ |x|}=\delta_{1i}\qquad \forall\
x\in B_1\setminus \{0\}.
$$
Let $\hat a_{11}=1$, $\hat a_{ii}=\epsilon$
for $2\le i\le n$, and $\hat a_{ij}=0$ for
$i\ne j$.
Define
$$
a_{ij}(x):= O(x)_{li} \hat a_{lm} O(x)_{mj}.
$$
Consider $u(x):=|x|^\alpha$, $v\equiv 0$, we have
$u>v$ in $B_1\setminus\{0\}$ and
\begin{eqnarray*}
F(x, \nabla^2 u)&:=& -a_{ij}(x)u_{ij}(x)
=-\alpha |x|^{\alpha-2}
a_{ij}(x)\left(
\delta_{ij}-(2-\alpha) \frac{x_ix_j}{ |x|^2}\right)\\
&=&- \alpha |x|^{\alpha-2}
\left( 1+(n-1)\epsilon-(2-\alpha)\right)
= 0=F(x, \nabla ^2 v).
\end{eqnarray*}
But (\ref{h}) holds.

\medskip

In Theorem \ref{thm1}, the singular set of $u$ is a point.
Our next theorem allows the singular set to be a
closed submanifold $E$ of dimension  $k\le n-2$.
A submanifold is called a closed submanifold if it is
a compact manifold without boundary.

\begin{thm}  For $n\ge 2$, let $F$ satisfy
(\ref{elliptic}), and   let $E\subset \Omega$ be a smooth closed
 submanifold of dimension $k$,  $0\le k\le n-2$.
  Assume that   $u\in
C^2(\Omega\setminus E)$ and $v\in C^2(\Omega)$ satisfy
\begin{equation}
\lambda_1(\nabla^2 u)+\cdots+\lambda_{k+2}(\nabla^2 u)\le 0,
\qquad \mbox{in}\  
\Omega\setminus E,
\label{ABC1}
\end{equation}
$$
u>v\qquad \mbox{in}\ \Omega\setminus E,
$$
\begin{equation}
F(x, u, \nabla u, \nabla^2 u)\ge F(x, v, \nabla v, \nabla^2v) \qquad
\mbox{in}\ \Omega\setminus E. \label{abcde} \end{equation} Then
\begin{equation}
\liminf_{ dist(x, E)\to 0}(u-v)(x)>0.
\label{aa2}
\end{equation}
\label{thm11new}
\end{thm}

\begin{rem} In the above theorem,
condition 
 (\ref{ABC1}) is only needed to be satisfied,
for some $\bar r>0$,
in
$E_{\bar r}\setminus E$,
$E_{\bar r}=\{x\ |\
dist(x, E)<\bar r\}$, since
we can apply the theorem with
$\Omega=E_{\bar r}$.
\end{rem}

Similar to Remark \ref{rembb1}, we have
\begin{rem}
If the assumption ``$u>v$ in
$\Omega\setminus E$'' in  the above theorem
is replaced by ``$u \ge v$ in
$\Omega\setminus E$'', then,
either $u\equiv v$ in $\Omega\setminus E$ or
(\ref{aa2}) holds.
\label{rembb2}
\end{rem}

Our proof of Theorem \ref{thm11new} makes use of the following
maximum principle for functions satisfying
 (\ref{ABC1}).

\begin{prop}
For $n\ge 2$, $-1\le k\le n-2$,  let $E$ be a smooth closed
$k-$dimensional manifold in $\Bbb R^n$ and
$\Omega\subset \Bbb R^n$ be a domain.
Assume that
$u\in C^2(\Omega\setminus E)\cap C^0(\overline \Omega\setminus E)$ satisfies
(\ref{ABC1}) and
$$
\inf_{\Omega\setminus E}u>-\infty.
$$
Then
$$
u\ge \inf_{ \partial \Omega \setminus E}u,\qquad \mbox{on}\ \Omega\setminus E.
$$
\label{prop1}
\end{prop}

Note that in the above, when $k=-1$, $E$ is understood as
$\emptyset$, the empty set;
 while for $k=0$, $E$ consists of finitely many points.

 Lemma \ref{lem0-1} is a special case of Proposition
\ref{prop1}, by taking $k=0$, $n\ge 2$,
$\Omega=B_1$ and $E=\{0\}$.

We introduce, for $k=1,2,\cdots, n$,
the notion of viscosity $k-$superaffine functions
which extends that  of superharmonic functions
and superaffine functions.

 We use  $LSC(\Omega)$
to denote the set of lower semi-continuous functions  in $\Omega$, i.e.
the set of those maps $u: \Omega\to \Bbb R\cup  \{-\infty\}$
satisfying
$$
\liminf_{x\to \bar x}u(x)\ge u(\bar x),\qquad \forall\ \bar x\in \Omega.
$$
Similarly,  we use  $USC(\Omega)$
to denote the set of upper semi-continuous functions in $\Omega$, i.e.
the set of those maps $u: \Omega\to \Bbb R\cup  \{\infty\}$
satisfying
$$
\limsup_{x\to \bar x}u(x)\le u(\bar x),\qquad \forall\ \bar x\in \Omega.
$$

\begin{definition}
Let $\Omega\subset \Bbb R^n$ be an open set, and let $u
\in LSC(\Omega)$ satisfy
\begin{equation}
\inf_\Omega u>-\infty.
\label{DD0}
\end{equation}

For $1\le k\le n$, we
say that $u$ is a viscosity $k-$superaffine function  in $\Omega$
 if:
$$
\varphi\in C^2(\Omega),\ \ \
\varphi \le u,\quad\mbox{in}\ \Omega,\qquad
\mbox{and}\
\varphi(\bar x)=u(\bar x),\ \mbox{for some}\
 \bar x\in \Omega
$$
imply
$$
\lambda_1(\nabla^2 \varphi)+ \lambda_2(\nabla^2 \varphi)
+
\cdots+ \lambda_k(\nabla^2  \varphi)\le 0,\quad
\mbox{at}\ \bar x.
$$
We say that $u$ is a $k-$subaffine function  in $\Omega$
if $-u$ is  a $k-$superaffine function  in $\Omega$.
\label{Definition 1.1}
\end{definition}
We use ${\cal A}_k(\Omega)$ to denote the set of viscosity 
$k-$superaffine functions in $\Omega$.

Clearly,  ${\cal A}_n(\Omega)$ is the set of superharmonic functions in
$\Omega$,  ${\cal A}_1(\Omega)$ is the set of superaffine
 functions in
$\Omega$, and
$$
{\cal A}_n(\Omega)\subset {\cal A}_{n-1}(\Omega)\subset
\cdots\subset {\cal A}_1(\Omega).
$$

In Appendix A, we include some results on 
viscosity $k-$superaffine functions, including
an extension of Proposition \ref{prop1} for
viscosity $k-$superaffine functions.

\medskip

Theorem \ref{thm11new} in the case $k=n-2$
can be extended  as
follows.

\begin{thm}
 For $n\ge 2$, let $F$ satisfy
(\ref{elliptic}), and   let $E\subset \Omega$
be a closed subset of zero (Newtonian) capacity.
Assume that    $u\in
LSC(\Omega\setminus E)$ and  $v\in C^2(\Omega)$ satisfy,
for some constant $C$,
\begin{equation}
\Delta u\le C \qquad \mbox{in}\ \Omega\setminus  E,
\label{abab0}
\end{equation}
\begin{equation}
u>v\qquad \mbox{in}\ \Omega\setminus E,
\label{13-14}
\end{equation}
and
\begin{equation}
F(x, u, \nabla u, \nabla^2 u)\ge F(x, v, \nabla v, \nabla^2v) \qquad
\mbox{in}\ \Omega\setminus E
\ \mbox{in the viscosity sense}.
\label{abab5}
\end{equation}
 Then
(\ref{aa2}) holds.
\label{thmnew1}
\end{thm}

\begin{rem}
Theorem \ref{thmnew1}
 gives complete answers to Question 1.1, 1.2 and 1.3
in \cite{L1}.
\end{rem}

\begin{rem}
In the above theorem, (\ref{13-14}) can be replaced by
$$
%\begin{equation}
u\ge v\ \mbox{in}\
\Omega\setminus E,
\ \mbox{and}\ u\ \mbox{is not
identically equal to}\ v\ 
\mbox{in}\ 
\Omega\setminus E,
$$
%\label{B1-new}
%\end{equation}
and
(\ref{abab0}) can be replaced by
\begin{equation}
a_{ij}(x)\partial_{ij}u+b_i(x)\partial_iu+c(x)u\le C,
\quad\mbox{in}\ \Omega\setminus E, \
\mbox{in viscosity sense}
\label{B1-newnew}
\end{equation}
where $(a_{ij}(x))$ is some H\"older continuous positive definite matrix
functions in $\Omega$, $b_i(x)$ and $c(x)$ are continuous
functions in $\Omega$, and $C$ is some constant;
see \cite{CLN2}.
\label{rem-B1}
\end{rem}

\begin{openproblem}
Extend Proposition \ref{prop1} and   Theorem
\ref{thm11new} for $1\le k< n-2$ to more general $E$.
\end{openproblem}

\subsection{}
As an application of  Theorem \ref{thm1},  we
establish some results on symmetry of positive
solutions  in a punctured ball
using  the method of moving planes.

For $p=(p_1, \cdots, p_n)\in  \Bbb R^n$ and
$M=(M_{ij})\in  {\cal S}^{ n\times n}$, we
use, in (\ref{cond1}) below,
 notation
$$
\hat p:= (-p_1, p_2,  \cdots, p_n),\qquad  \widehat M:=(\widehat M_{ij})
$$
with
$
\widehat M_{1j}=-M_{1j}
$ for $2\le j\le n$;
$M_{i1}=-M_{i1}$ for $2\le i\le n$;
$\widehat M_{11}=M_{11}$;
$\widehat M_{ij}=M_{ij}$
for $2\le i, j\le n$.

\begin{thm} For $n\ge 1$,
let $F\in C^1(\overline B_1\setminus\{0\}\times \Bbb R\times  \Bbb R^n\times
 {\cal S}^{ n\times n})$ satisfy 
\begin{equation}
\left( \frac{\partial F}{\partial M_{ij} }(x, s, p, M)\right)
>0\qquad \forall\ (x, s, p, M)\in
\overline B_1\setminus \{0\}
 \times \Bbb R\times \Bbb R^n\times {\cal S}^{n\times n},
\label{Q1-1}
\end{equation}
(the reader is warned that 
$-F$ satisfies (\ref{elliptic})).
For all
 $x\in B_1$,  $0<x_1<1$, $-1<\hat x_1\le
x_1$, $x_1+\hat x_1>0$,  $p\in   \Bbb R^n$,
$M\in  {\cal S}^{ n\times n}$,
\begin{equation}
F(\hat x_1, x_2, \cdots, x_n, s, \hat p, \widehat M)
\ge
F(x_1,  x_2, \cdots, x_n, s, p, M).
\label{cond1}
\end{equation}
Assume that $u\in C^2(\overline B_1\setminus\{0\})$
satisfies  (\ref{superaffine}) for some $\bar r>0$,
\begin{equation}
F(x, u, \nabla u, \nabla ^2 u)=0,\ \ u>0,\qquad \mbox{in}\ B_1\setminus \{0\}
\label{cond2}
\end{equation}
and
$$
u=0\qquad \mbox{on}\ \partial B_1.
$$
Then
\begin{equation}
\frac { \partial u}{ \partial x_1}(x)<0\qquad
\forall\ x\in B_1, \ 0<x_1<1.
\label{mono}
\end{equation}
 \label{thm3}
\end{thm}

The theorem is proved in section 4.  An immediate consequence is
\begin{cor} \label{cor-Q2}
If in addition, $F$ is symmetric under reflection
in the plane $\{x_1=0\}$, then
$u$ is symmetric in $x_1$.
\end{cor}

Another consequence, after applying the theorem in every direction,
is

\begin{cor} For $n\ge 1$,
let $F(|x|, s, p, M)$ satisfy (\ref{elliptic}),
and for some $\epsilon>0$,
$$
\frac {\partial }{\partial r}F(r, s, p, M)\ge 0\qquad
\forall\ 0<r<1+\epsilon,  s>0, p\in \Bbb R^n, M\in  {\cal S}^{ n\times n},
$$
 and
$$
F(r, s, Op, O^t MO)\equiv F(r, s, p, M),\quad \forall\
0<r<1+\epsilon,
s>0, p\in \Bbb R^n, M\in {\cal S}^{ n\times n}, O\in O(n),
$$
 and let
$u\in C^2(\overline B_1\setminus\{0\})$
satisfy (\ref{superaffine})
for some $\bar r>0$,
$$
F(|x|, u, \nabla u, \nabla ^2 u)=0,\ \ u>0,\qquad \mbox{in}\ B_1\setminus \{0\}
$$
and
$$
u=0\qquad \mbox{on}\ \partial B_1.
$$
Then $u$ is radially symmetric about the origin and $u'(r)<0$ for
$0<r<1$. \label{cor3}
\end{cor}

Condition (\ref{superaffine})
cannot  be dropped in Corollary \ref{cor3}, as shown by the
following examples.

\begin{exam}
Let
$u\in C^2([-1, 1]\setminus\{0\})$ be defined
by $u(x)=1+x$ if $-1\le x\le 0$,
$u(x)=\frac 12 -\frac x2$ if $0<x\le 1$.
Then,
$$
u''=0, u>0 \quad\mbox{in}\ (-1, 1)\setminus\{0\},
$$
but $u$ is not symmetric about the origin.
This shows that condition (\ref{superaffine})
cannot be dropped in Corollary \ref{cor3} in dimension $n=1$.
\label{example4}
\end{exam}

Here is another example 
 in  dimension $n=1$.
\begin{exam}
Let
$$
u(x)=\sin(\pi |x|),\qquad 0<|x|<1.
$$
Then
$$
u''+\pi^2 u=0, \ u>0, \qquad \mbox{in}\ (-1, 1)\setminus\{0\},
$$
but $u$ is not monotone in $(0, 1)$.
\label{example5}
\end{exam}

In a lecture about this paper
in November 2008, the last author
said that for $n>1$ he did not know if the condition
(\ref{superaffine})
could be dropped.  In particular, he asked if one had spherical symmetry
and monotonicity for a solution of
$$
\left\{
\begin{array}{rl}
\Delta u+f(u)&=0, \ \ \ u>0,\qquad B_1\setminus\{0\},\\
u&=0,  \qquad \qquad \qquad \partial B_1,
\end{array}
\right.
$$
$f$ locally Lipschitz --- with no
further assumption.  Immediately after the talk, Susanna Terracini
said she could prove it and she showed him the proof.
It made use of an idea she had used in \cite{T}
which treats equations with singular potential in all of space.  Her
proof uses a nice variant of the argument of
\cite{BN} and works also if $f$ is
permited
to depend also on $x$,
$f=f(x,u)$, satisfying 
(\ref{cond1}).

In a later paper \cite{CLN2},  we will make use of her variant
to prove symmetry and monotonicity
in problems of the form
(\ref{cond2}),  but under other conditions on $F$.

Recently we found an example for dimension $n\ge 2$
showing that condition (\ref{superaffine})
may not be dropped; Consider
$$
{\cal M}^+(M):=\sum_{i=1}^n h(\lambda_i(M)),\qquad
M\in {\cal S}^{n\times n},
$$
where $\{\lambda_i(M)\}$ denotes eigenvalues of
$M$ and
$$
h(s):=
\left\{
\begin{array}{ll}
1, & s\ge 0,\\
\frac{n-1}{1-\alpha}, & s<0.
\end{array}
\right.
$$

${\cal M}^+(\nabla ^2 u)$, a Pucci extremal operator,
is uniformly elliptic and is concave and Lipschitz
 in $\nabla^2 u$.
See e.g. section 2.2 of \cite{CC}
for the  Pucci extremal operators and
their properties.

In the following we will consider radially symmetric functions $u(x)$, and
we will slightly abuse notations by writing
$u(x)=u(r)$, $r=|x|$.

For $0<\alpha<1$, let
$$
u(r)=
\left\{
\begin{array}{ll}
r^\alpha-\frac \alpha 2 r^2, & 0<r\le 1,\\
-\frac{ \alpha (2-\alpha) }n\left( \frac 12 r^2+\frac 1{n-2} r^{2-n}\right)+
\frac {(2-\alpha) (n-2+\alpha) }{ 2(n-2) }, & r>1.
\end{array}
\right.
$$
A calculation gives
$$
\frac{u'(r)}r=
\left\{
\begin{array}{ll}
\alpha r^{\alpha-2}- \alpha , & 0<r\le 1,\\
- \frac{ \alpha (2-\alpha) }n (1-r^{-n}),  & r>1,
\end{array}
\right.
$$
and
$$
u''(r)=
\left\{
\begin{array}{ll}
-\alpha(1-\alpha)r^{\alpha-2}
-\alpha, & 0<r\le 1,\\
- \frac{ \alpha (2-\alpha) }n \left( 1+(n-1) r^{-n}\right),&
r>1.
\end{array}
\right.
$$
The function $u$ is in $C^{2,1}(0, \infty)$ and satisfies,
for some $\bar r>1$,
$$
u>0\ \mbox{in}\ (0, \bar r), \qquad u(0)=u(\bar r)=0.
$$

 The eigenvalues
of the Hessian $\nabla ^2 u$ are
$$
\lambda_1(\nabla^2 u)= u''(r),
\qquad \lambda_2(\nabla^2 u)=\cdots= \lambda_n(\nabla^2 u)=
\frac{ u'(r)}r.
$$
Since
$
u''<0$ in $(0, \infty)$, while
$u'>0$ in $(0, 1)$ and $u'<0$ in $(1, \infty)$,
we have
\begin{eqnarray*}
{\cal M}^+(\nabla^2 u)(x)&=&
\frac{n-1} {1-\alpha} \lambda_1(\nabla^2 u)
+\sum_{i=2}^n  \lambda_i(\nabla^2 u)=
\frac{n-1} {1-\alpha} u''(r)+
(n-1) \frac{ u'(r)}r\\
&=&
-\frac{  (n-1) \alpha (2-\alpha) }{1-\alpha},\qquad
\qquad 0<|x|\le 1,
\end{eqnarray*}
and
\begin{eqnarray*}
{\cal M}^+(\nabla^2 u)(x)&=&
\frac{n-1} {1-\alpha}
 \sum_{i=1}^n  \lambda_i(\nabla^2 u)=
\frac{n-1} {1-\alpha} \left( u''(r)+\frac {n-1}r u'(r)\right)\\
&=&
-\frac{  (n-1) \alpha (2-\alpha) }{1-\alpha},\qquad
\qquad |x|>1.
\end{eqnarray*}

\begin{exam} For $n\ge 2$, $0<\alpha<1$, let $u$ and ${\cal M}^+$
be as above.  Then
$u\in C^{2,1}(\overline{B_{\bar r}(0)}\setminus\{0\})$ satisfies
$$
{\cal M}^+(\nabla^2 u)=
-\frac{  (n-1) \alpha (2-\alpha) }{1-\alpha},\qquad
\mbox{in}\ B_{\bar r}(0) \setminus\{0\},
$$
and
$$
u>0\ \  \mbox{in}\ B_{\bar r}(0) \setminus\{0\},
\qquad u=0\ \  \mbox{on}\ \partial B_{\bar r}(0).
$$
But $u'(r)$ changes signs in $(0, \bar r)$.
\end{exam}

\subsection{More general equations}
We  extend Theorem \ref{thm11new}
to more general operators which are motivated by
conformally invariant operators.
An operator $F(\cdot, u, \nabla u, \nabla^2u)$  is conformally  invariant
if
$$
F(\cdot, u_\psi, \nabla u_\psi, \nabla^2 u_\psi)\equiv
F(\cdot, u, \nabla u, \nabla^2 u)\circ \psi\qquad\mbox{on}\
\Bbb R^n
$$
for all  M\"obius transformation
$\psi$ and all  positive functions  $u\in C^2(\Bbb R^n)$, where
$
u_\psi:=|J_\psi|^{ \frac {n-2}{2n} }(u\circ \psi)$
and  $J_\psi$ denotes the Jacobian of $\psi$.
It was proved in \cite{LL1},
see Theorem 1.2 there,  that conformally invariant
operators are all
of the form $f(\lambda(A^u))$, where $f(\lambda)$ is a symmetric function of
$\lambda=(\lambda_1, \cdots, \lambda_n)$, and
$\lambda(A^u)$ denotes the eigenvalues of the
$n\times n$ symmetric matrix function
\begin{equation}
A^u:= -\frac{2}{n-2}u^{  -\frac {n+2}{n-2} }
\nabla^2u+ \frac{2n}{(n-2)^2}u^ { -\frac {2n}{n-2} }
\nabla u\otimes\nabla u-\frac{2}{(n-2)^2} u^ { -\frac {2n}{n-2} }
|\nabla u|^2I,
\label{abc1}
\end{equation}
and  $I$ is the $n\times n$ identity matrix.
After setting
$$
w=u^{ -\frac 2{n-2}},
$$
we have
$$
%\begin{equation}
A^u=A_w:= w \nabla^2 w- \frac {|\nabla w|^2}{2}I,
$$
%\end{equation}
which has been called in  some literature the conformal Hessian of $w$.
The ellipticity of $f(\lambda(A^u))$ amounts to
$\displaystyle{
\frac{\partial f}{\partial \lambda_i}>0, 1\le i\le n}$,
in the relevant region.

If $F(x, u, \nabla u, \nabla^2 u)$ is a conformally invariant
second order elliptic
operator,
(\ref{aa1}) was proved in \cite{L1} if $u$ satisfies, instead of
(\ref{superaffine}), a stronger assumption $\Delta u\le 0$
in $B_3\setminus \{0\}$; see Theorem 1.6 there.
As explained in the introduction of \cite{LL2},
the general Liouville theorem  for
 conformally invariant  elliptic
operators of second order  established there
(Theorem 1.3 in the paper)  would have followed from the
proof of Theorem 1.4 in \cite{LL1} if
the above mentioned theorem in  \cite{L1}
had been known at the time.
The proof of this theorem in \cite{L1}
makes use of ideas in the proof of the Liouville theorem in
\cite{LL2}.  Given (\ref{aa1})  for
conformally invariant elliptic
operators, one can also adapt the moving
planes procedure in \cite{GNN} and \cite{CGS}
to show that  all solutions are
radially symmetric about some point, and then
classify  radially symmetric entire solutions.

In  \cite{L1}, (\ref{aa1}) was also established
for some special classes of $F(x, u, \nabla u, \nabla^2 u)$ which include
$F(\nabla u, \nabla^2 u)$, again for superharmonic $u$;
see Theorem 1.7 and Corollary 1.6 in \cite{L1}.
Our proof of Theorem \ref{thm1},  very different from
the arguments in  \cite{L1},
is also extended to  the situations
where the singular set $\{0\}$ is replaced
by a  closed submanifold
of dimension $k\le n-2$
(see Theorem \ref{thm11new}).
Our proof of Theorem
\ref{thmnew1}, different from
that of Theorem \ref{thm1} and the arguments  in  \cite{L1},
makes use of the Alexandrov-Bakelman-Pucci estimate;
see e.g. section 3 of \cite{CC} and section 9.1 of
\cite{GT} for the ABP estimate.

Here we consider
\begin{equation}
\widetilde A^u:= -\nabla ^2 u +\widetilde L(\cdot, u, \nabla u)
\label{3.1}
\end{equation}
where $\widetilde L$
satisfies
\begin{equation}
\widetilde L \in C^{0, 1}_{loc}(\Bbb R^n\times \Bbb R\times \Bbb
R^n) .\label{3.2}
\end{equation}

One example of such $\widetilde L$
is,
$$
\widetilde L(\cdot, u, \nabla u)=
\frac n{(n-2)u}  \nabla u\otimes\nabla u
-\frac 1{(n-2)u}
|\nabla u|^2 I.
$$
With  this $\widetilde L$, $\widetilde A^u$ takes the form
\begin{equation}
\widetilde A^u= \frac{n-2} 2 u^{ \frac {n+2}{n-2} }
A^u,
\label{abc2}
\end{equation}
where $A^u$ is the conformally
invariant one given in (\ref{abc1}).

Fully nonlinear elliptic equations of second order
with the Hessian $\nabla ^2 u$ in appropriate regions
of ${\cal S}^{n\times n}$ have been investigated in  
the classical paper \cite{CNS} and many subsequent ones.

Let
   $U\subset {\cal S}^{n\times n}$
   be  an open set,
we assume that
\begin{equation} M+N\in U,\qquad \forall\
M\in U, N\in  {\cal S}^{n\times n}_+.
\label{YY1}
\end{equation}

For example, for $1\le k\le n$, let
$$
\sigma_k(\lambda)=\sum_{1\le i_1<\cdots <i_k\le n}\lambda_{i_1}
\cdots \lambda_{i_k}
$$
be the $k-$th elementary  symmetric function and let
$
\Gamma_k$ be the connected component
of $\{\lambda\in \Bbb R^n\ |\
\sigma_k(\lambda)>0\}$ containing the positive cone $\Gamma_n$.
Then
\begin{equation}
U_k:=\{ M\in {\cal S}^{ n\times n}\ |\
\lambda(M)\in \Gamma_k\}
\label{13-aa}
\end{equation}
satisfies (\ref{YY1}) ---
see e.g. \cite{CNS}.

\begin{thm}
For $n\ge 2$, let $\Omega$ be a domain
 in $\Bbb R^n$ containing the origin
$\{0\}$, 
 $U$ be an open subset of
$  {\cal S}^{n\times n} $
satisfying  (\ref{YY1}),
 and let
$\widetilde A^u$ satisfy (\ref{3.1}) and (\ref{3.2}).  Assume that
$u\in C^2(\Omega\setminus \{0\})$  satisfies (\ref{superaffine})
for some $\bar r>0$,
    $v\in C^2(\Omega)$,
$$
%\begin{equation}
\widetilde A^u\in
\overline U\qquad \mbox{in}\
\Omega\setminus \{0\},
$$
%\label{2-1new}
%\end{equation}
$$
%\begin{equation}
\widetilde A^v\in  {\cal S}^{n\times n} \setminus U\qquad \mbox{in}\
\Omega,
$$
% \label{2-2new}
%\end{equation}
\begin{equation}
u>v\qquad \mbox{in}\
\Omega\setminus \{0\}.
\label{2-3new}
\end{equation}
Then
\begin{equation}
\liminf_{ |x|\to 0}(u-v)(x)>0.
\label{3-1new}
\end{equation}
\label{thm2-1new}
\end{thm}

\begin{question}
Does (\ref{3-1new}) still hold if
(\ref{2-3new}) in Theorem \ref{thm2-1new}
is replaced  by
 a weaker assumption
\begin{equation}
u\ge v\ \mbox{in}\ \Omega\setminus \{0\},\ \ \mbox{and}\
 \liminf_{dist(x, \partial \Omega)\to 0}
(u-v)(x)>0?
\label{27new}
\end{equation}
\label{bbbnew}
\end{question}

For  an open subset  $U$  of
$  {\cal S}^{n\times n} $, we introduce, for  $1\le k\le n-1$,

\noindent {\it Property P$_{k}$}:
Given any  constant $C>0$, there exist some
constants $\bar \delta=\bar \delta(C)>0,
\bar \epsilon=\bar \epsilon(C)>0$ such that
for any $0<\epsilon<\bar \epsilon$, and 
$
 N\in {\cal S}^{n\times n}, M\in \overline U,  M\le CI,
\lambda_i(N)=-\bar \delta, 1\le i\le k,
\lambda_j(N)=1, k+1\le j\le n
 \
\mbox{we have}\ M+\epsilon N\in U.
$

\begin{rem} 
Property P$_k$ is satisfied by an open convex  subset $U$ of
${\cal S}^{n\times n}$ if it further
satisfies one of the following:

\medskip

\noindent (i)\  $\partial U\cap
\{M\ |\ M\le CI\}$ is compact for any constant $C$,
and, for any $M\in \partial U$, there exists
a constant $a>0$ such that
$M+N\in U$ if 
eigenvalues of $N$ 
consist of $k$ zeros and $n-k$ $a^{'}$s.

\medskip

\noindent (ii)\ 
$M\in U$ implies $aM\in U$ for any positive 
constant $a$, 
and $N$ with  $k$ zeros and $(n-k)$ ones as
eigenvalues belongs to $U$.

\medskip

\noindent (iii)
$U=U_l$, as in (\ref{13-aa}),  for some $1\le l\le n-k$.
\label{rem-13}
\end{rem}
See the end of section 6 for an explanation of the above 
remark.

In Theorem \ref{thm2-1new}, the singular set is a point $\{0\}$.
The following theorem allows singular sets to be 
 submanifolds of dimension $k$,  $1\le k\le n-2$.

\begin{thm}
For $n\ge 2$, let $\Omega\subset \Bbb R^n$ be a domain,
  $E\subset \Omega$ be a smooth closed
 submanifold of dimension $k$,  $1\le k\le n-2$,
 $U$ be an open subset of
$  {\cal S}^{n\times n} $
satisfying  (\ref{YY1})
and Property P$_{k}$,
 and let
$\widetilde A^u$ satisfy (\ref{3.1}) and (\ref{3.2}).  Assume that
$u\in C^2(\Omega\setminus E)$ satisfies
 (\ref{ABC1}),   $v\in C^2(\Omega)$, and
for some $\bar r>0$,
\begin{equation}
\widetilde A^u\in
\overline U\qquad \mbox{in}\
\Omega\setminus E,
\label{2-1}
\end{equation}
\begin{equation}
\widetilde A^v\in  {\cal S}^{n\times n} \setminus U\qquad \mbox{in}\
\Omega, \label{2-2}
\end{equation}
\begin{equation}
u>v\qquad \mbox{in}\
\Omega\setminus E.
\label{2-3}
\end{equation}
Then
\begin{equation}
\liminf_{ dist(x, E)\to 0}(u-v)(x)>0.
\label{3-1}
\end{equation}
\label{thm2-1}
\end{thm}

Similar to Question \ref{bbbnew}, we have
\begin{question}
Does (\ref{3-1}) still hold if
(\ref{2-3}) in Theorem \ref{thm2-1}
is replaced  by
 a weaker assumption
$$
u\ge v\ \mbox{in}\ \Omega\setminus E,\ \ \mbox{and}\
 \liminf_{dist(x, \partial \Omega)\to 0}
(u-v)(x)>0?
$$
\label{bbb}
\end{question}

We suspect that the answers to Question \ref{bbbnew}
and \ref{bbb} are ``No'', and some structure assumptions
on $\widetilde L$ are needed in order to have
an affirmative answer.
These issues will be addressed in a subsequent paper \cite{CLN1}.

If $\widetilde A^u$ is the conformally invariant one given by (\ref{abc2}),
and if condition 
(\ref{superaffine})
 is replaced by a stronger
assumption ``$\Delta u\le 0$ in $B_{\bar r}\setminus\{0\}$
 for some $\bar r>0$'',
then the answer to Question \ref{bbbnew} is ``Yes'', as
proved in
\cite{L4}.
  The assumption (\ref{27new}) is indeed weaker than (\ref{2-3new}) since
the strong maximum principle
might not hold --- indeed, as shown in \cite{LN},
both  the
strong maximum principle and the Hopf Lemma fail
for $U=U_k$, $2\le k\le n$:  there exist $u_1\in C^\infty(B_1)$ such that
$u_1$ and $u_2\equiv 1$ satisfy
$0<u_1<u_2$ in $B_1\setminus \{x_1=0\}$, $u_1=u_2$
on $\{x_1=0\}$,  $A^{u_i}\in \partial U_k$ for $i=1,2$,
but $\displaystyle{
\partial_{x_1}(u_1-u_2)(0)=0}$.

The above mentioned result in \cite{L4}
 was established under a weaker
$C^{1,1}$ regularity assumption on $u$ and $v$.
The regularity assumption was further
weakened in \cite{L2}
 to locally Lipschitz  viscosity solutions
of
(\ref{2-1}) and (\ref{2-2}), which is crucial
in the proof there of the local gradient
estimates to general conformally invariant second order elliptic
equations without concavity assumptions on the equations.

The proofs of the above mentioned results in
 \cite{L4} and  \cite{L2}, as well as
that of the previously mentioned results
in \cite{L1} for conformally invariant
elliptic operators, make use of the following

  \medskip

\noindent{\bf Lemma A}\
{\it (\cite{LL2})\
For $n\ge 2$, $B_1\subset\Bbb {R}^n$, let $u\in L^1_{loc}(B_1\setminus
\{0\})$ satisfy
$$
\liminf_{x\to 0} u(x)=0,
$$
\[
\Delta u\le 0\qquad\mbox{in}~B_1\setminus \{0\}
\]
in the distribution sense.
Let $P$  denote the set of vectors
 $p$ in $\Bbb R^n$ satisfying
\begin{equation}
u(x)\ge p\cdot x+
\circ(|x|)\qquad\mbox{in}\
B_1\setminus \{0\}.
\label{(3)}
\end{equation}
Then $P$ contains at most one element.
}

In the above, $\Delta u\le 0$ can be replaced 
by $\Delta u\le C$ for some constant $C$ --- 
since we can work with $u(x)-\frac C{2n} |x|^2$ instead of $u$.

In Appendix F we give an improvement of the above lemma, see
Lemma \ref{lem-F3} there, which is used in the
proof of Theorem \ref{thmnew1} (see the proof of
Lemma \ref{lem-35}).
We also mention an extension of the  lemma for functions satisfying
(\ref{superaffine}) instead of superharmonic
functions in the punctured ball.
It is not used in this paper.

\begin{lem}
Let $u$ be defined in $B_1\setminus \{0\}$
in $\Bbb R^n$, $n\ge 2$, with
$$
%\begin{equation}
\liminf_{x\to 0}u(x)=0.
$$
%\label{(1)}
%\end{equation}
Suppose that $u$ satisfies (\ref{superaffine})
with $\bar r=1$.
Let $P$  denote the set of vectors
 $p$ in $\Bbb R^n$ satisfying
(\ref{(3)}).
Then there are at most $n$ vectors
$p_0, \cdots, p_{n-1}$ in
$P$ such that
\begin{equation}
p_1-p_0, \cdots, p_{n-1}-p_0
\ \mbox{are linearly independent},
\label{(4)}
\end{equation}
i.e. the convex hull of all the vectors in $P$
is at most $(n-1)$ dimension.
\label{abc}
\end{lem}

\begin{rem}
Condition (\ref{(4)}) is equivalent to:
for any $i\le n-1$,
$\{p_j-p_i\ |\ 0\le j\le n-1, j\ne i\}$
are linearly independent.
\end{rem}

The paper is organized as follows.
Theorem \ref{thm1} and
Theorem \ref{thm11new}
are proved in section 2,
Theorem \ref{thmnew1}
is proved in section 3,
Theorem \ref{thm3} is proved in section 4,
Theorem \ref{thm2-1new} is proved in section 5, and
 Theorem \ref{thm2-1}  is proved in section 6 --- the proof is
technical and it is suggested to omit
this section on a first reading.
Appendix A contains some properties of distance functions,
Appendix B contains some linear algebra lemmas.
Proposition \ref{prop1} and an extension of it,
Theorem  \ref{Theorem 1.1},  are proved in Appendix C and Appendix D,
 Lemma \ref{abc} is proved in Appendix E.
Appendix F contains some improvement of Lemma A.

\section{Proof of   Theorem \ref{thm1} and
Theorem \ref{thm11new}}

\subsection{}
\noindent{\bf Proof of Theorem \ref{thm1}.}\
 We prove it by
contradiction.  Suppose the contrary, then
\begin{equation}
\liminf_{x\to 0}
(u-v)(x)=0. \label{zeronew} \end{equation}

Let
$$
v_\epsilon(x):= v(x)+\epsilon |x|.
$$

\begin{lem}
For any $\epsilon>0$,
there exists $0<\bar r_\epsilon<1$ such that
$\displaystyle{
\inf_{ \partial B_r  }(u-v_\epsilon)<0
}
$ for
all $0<r<\bar r_\epsilon$.
\label{lemma0}
\end{lem}

\noindent{\bf Proof.}\
By  condition
 (\ref{superaffine})  and   (\ref{zeronew}),
\begin{eqnarray*}
\inf_{ x\in  \partial  B_r }
[u(x)-(v(0)+\nabla v(0)\cdot x)]
&=&\inf_{  B_r \setminus\{0\} }[u(x)-(v(0)+\nabla v(0)\cdot x)]
\\
&\le& \liminf_{ x\to 0}[ u(x)-(v(0)+\nabla v(0)\cdot x) ]
=0.
\end{eqnarray*}
Thus
\begin{eqnarray*}
\inf_{ \partial B_r }(u-v_\epsilon)
&\le &
\inf_{ x\in  \partial  B_r }
[u(x)-(v(0)+\nabla v(0)\cdot x)]
+\sup_{ x\in  \partial  B_r }
[(v(0)+\nabla v(0)\cdot x) -v(x)
-\epsilon |x|]
\\
&\le &
\frac 12 (\sup_{B_1}\|\nabla^2 v\|)r^2 -\epsilon r.
\end{eqnarray*}
Lemma \ref{lemma0} follows.

\bigskip

By Lemma \ref{lemma0},
$$
\lambda(\epsilon):=-\inf_{ B_3\setminus \{0\} }(u-v_\epsilon)>0.
$$
 Since
$$
 \liminf_{ x\to 0}
 (u-v_\epsilon)(x)= \liminf_{ x\to 0}
 (u-v)(x)= 0,
 $$
 there
exists $x_\epsilon\in \overline B_3\setminus \{0\}$ such that
\begin{equation}
u-\tilde v_\epsilon\ge 0\ \mbox{in}\
B_3\setminus \{0\}, \qquad
(u-\tilde v_\epsilon)(x_\epsilon)=0,
\label{13-1new}
\end{equation}
where
 $$\tilde v_\epsilon(x):=v_\epsilon(x)-\lambda(\epsilon)
=v(x)+\epsilon |x|-\lambda(\epsilon).
$$
Using the positivity of $u-v$ in $\overline B_3\setminus\{0\}$,
we obtain from the above that
\begin{equation}
\lambda(\epsilon)=-(u-v_\epsilon)(x_\epsilon)=
v(x_\epsilon)-u(x_\epsilon) +\epsilon |x_\epsilon|\le \epsilon
|x_\epsilon|,\label{12-1}
\end{equation}
and
$$
\lim_{\epsilon\to 0}
|x_\epsilon|=0.
$$

We will show, for small $\epsilon>0$,  that
\begin{equation}
F(x_\epsilon, v(x_\epsilon), \nabla v(x_\epsilon), \nabla
^2v(x_\epsilon))
> F(x_\epsilon, \tilde v_\epsilon(x_\epsilon),
\nabla \tilde v_\epsilon(x_\epsilon),
\nabla^2 \tilde v_\epsilon(x_\epsilon)).  \label{13-3new}
\end{equation}
This would lead to a contradiction.  Indeed this would
 imply, together
with (\ref{abcde}), that
$$
F(\cdot, u, \nabla u, \nabla^2 u)(x_\epsilon)
> F(\cdot, \tilde v_\epsilon, \nabla \tilde
v_\epsilon, \nabla ^2 \tilde v_\epsilon)(x_\epsilon).
$$
We deduce from (\ref{13-1new}) that
 $\nabla
u(x_\epsilon)=\nabla \tilde v_\epsilon(x_\epsilon)$, $(\nabla^2
u(x_\epsilon))\ge (\nabla ^2\tilde v_\epsilon(x_\epsilon))$, and,
using the ellipticity assumption (\ref{elliptic}),
$$
F(\cdot, u, \nabla u, \nabla^2 u)(x_\epsilon)\le F(\cdot, \tilde v_\epsilon,
\nabla \tilde v_\epsilon, \nabla ^2 \tilde v_\epsilon)(x_\epsilon).
$$
A contradiction.

Now we prove (\ref{13-3new}) for small $\epsilon$.
By (\ref{12-1}),
\begin{equation}
|\tilde v_\epsilon-v|+
 |\nabla (\tilde
v_\epsilon-v)|=O(\epsilon), \label{15-1new}
\end{equation}
and therefore
\begin{eqnarray}
F(\cdot, v, \nabla v, \nabla ^2v)(x_\epsilon)
 &=&F(\cdot, \tilde v_\epsilon+O(\epsilon),
\nabla \tilde
v_\epsilon+O(\epsilon), \nabla ^2 v)(x_\epsilon) \nonumber\\
&\ge &F(\cdot, \tilde v_\epsilon, \nabla \tilde v_\epsilon, \nabla ^2
v)(x_\epsilon)-C\epsilon,
\label{16-1new}
\end{eqnarray}
where $C$ is independent of $\epsilon$.

By the mean value theorem,
\begin{eqnarray*}
&&
F(\cdot,  \tilde v_\epsilon, \nabla \tilde v_\epsilon, \nabla ^2
v)(x_\epsilon)-
 F(\cdot,  \tilde v_\epsilon, \nabla \tilde v_\epsilon, \nabla ^2
\tilde v_\epsilon)(x_\epsilon)
\\
&=& \int_0^1 F_{ij}(x_\epsilon, \tilde v_\epsilon(x_\epsilon),
 \nabla \tilde v_\epsilon(x_\epsilon),
\nabla ^2
v(x_\epsilon)+t[ \nabla ^2
\tilde v_\epsilon(x_\epsilon)
-\nabla ^2
v(x_\epsilon)])dt (\partial_{ij} v(x_\epsilon)
-\partial_{ij} \tilde v_\epsilon(x_\epsilon)).
\end{eqnarray*}

Consider the compact set
$$
S:=\{ (x,s,p,M)\ |\ |x|\le 1, \|(s, p, M)-(v(x), \nabla v(x),
\nabla ^2 v(x))\|\le 1\}.
$$
By  the ellipticity assumption (\ref{elliptic}),
there exists some constant $b>0$ such that
$$
\left(-\frac{\partial F}{\partial M_{ij} }\right )\ge bI
\qquad\mbox{on}\ S,
$$
where $I$ is the $n\times n$ identity matrix.

A calculation gives, for some  $O(x)\in O(n)$,
\begin{equation}
\nabla^2 \tilde v_\epsilon(x)-\nabla ^2 v(x)=
\epsilon \nabla ^2 (|x|) =
\epsilon O(x)^t
\left(diag(0, \frac 1 {|x|},  \cdots, \frac 1 {|x|})
\right) O(x).
\label{CC1}
\end{equation}

Define
$$
\bar t_\epsilon:=
\left\{
\begin{array}{rl}
1, & \mbox{if}\  \|\nabla^2 \tilde v_\epsilon(x_\epsilon)-
\nabla ^2 v(x_\epsilon)\|\le  1/\sqrt{3},
\\
1/ (\sqrt{3}  \|\nabla^2 \tilde v_\epsilon(x_\epsilon)-
\nabla ^2 v(x_\epsilon)\|),
 & \mbox{if}\ \|\nabla^2 \tilde v_\epsilon(x_\epsilon)-
\nabla ^2 v(x_\epsilon)\|>1/\sqrt{3}.
\end{array}
\right.
$$

In view of (\ref{15-1new}) and the definition of
$\bar t_\epsilon$,
$$
(x_\epsilon, \tilde v_\epsilon(x_\epsilon),
 \nabla \tilde v_\epsilon(x_\epsilon),
\nabla ^2
v(x_\epsilon)+t[ \nabla ^2
\tilde v_\epsilon(x_\epsilon)
-\nabla ^2
v(x_\epsilon)])\in S,\qquad
\forall\ 0\le t\le \bar t_\epsilon.
$$
It follows that
\begin{eqnarray*}
&&
F(\cdot,  \tilde v_\epsilon, \nabla \tilde v_\epsilon, \nabla ^2
v)(x_\epsilon)-
 F(\cdot,  \tilde v_\epsilon, \nabla \tilde v_\epsilon, \nabla ^2
\tilde v_\epsilon)(x_\epsilon)
\\
&\ge &
 \int_0^{\bar t_\epsilon} F_{ij}(x_\epsilon, \tilde v_\epsilon(x_\epsilon),
 \nabla \tilde v_\epsilon(x_\epsilon),
\nabla ^2
v(x_\epsilon)+t[ \nabla ^2
\tilde v_\epsilon(x_\epsilon)
-\nabla ^2
v(x_\epsilon)])dt (\partial_{ij} v(x_\epsilon)
-\partial_{ij} \tilde v_\epsilon(x_\epsilon))\\
&\ge &
b \bar t_\epsilon \sum_{i=1}^n
(\partial_{ii}\tilde  v_\epsilon(x_\epsilon) -\partial_{ii} v(x_\epsilon))
\ge
b \bar t_\epsilon \|\nabla ^2\tilde   v_\epsilon(x_\epsilon) -
\nabla ^2  v(x_\epsilon)\|\\
&\ge &
b\min \{  \|\nabla ^2 \tilde  v_\epsilon(x_\epsilon) -
\nabla ^2  v(x_\epsilon)\|, 1/\sqrt{3}\}\ge
b \min\{\frac \epsilon{|x_\epsilon|}, 1/\sqrt{3}\}.
\end{eqnarray*}
Since $|x_\epsilon|\to 0$, estimate
(\ref{13-3new}) follows from
  (\ref{16-1new}) and
the above.
 Theorem \ref{thm1} is established.

\subsection{}
Before proving Theorem \ref{thm11new} we first take up
a simpler case:

\noindent{\bf Proof of Theorem \ref{thm11new} for 
k=
n-2.}\
Shrinking $\Omega$ slightly, we may assume that
$u>v$
 on $\overline \Omega\setminus E$,  $u\in C^2(\overline \Omega\setminus E)$,
and $v\in  C^2(\overline \Omega)$.
 We prove it by
contradiction.   Suppose the contrary, then for some $\bar
x\in E$,
\begin{equation}
\liminf_{x\in \Omega\setminus E,
 x\to \bar x}(u-v)(x)=0. \label{zeronewnew} \end{equation}
Without loss of generality, $\bar x=0$.

For $C_1:=\frac 1{2n} \sup _{\Omega}(-\Delta v)$.
$$
\Delta (u-v-C_1|x|^2)\le 0,\qquad x\in \Omega\setminus E.
$$
Since $Cap(E)$, the capacity of $E$, is equal to $0$ and
$u-v-C_1|x|^2$ is bounded from below on
$\Omega\setminus E$,
$$
\Delta (u-v-C_1|x|^2)\le 0,\qquad \mbox{in}\ \Omega\ \mbox{in
the distribution sense}.
$$
It follows,  in view of (\ref{zeronewnew}),  that
\begin{equation}
\int_{ \partial B_r(0)}(u-v-C_1|x|^2)\le 0,\qquad 0<r<1.
\label{C5-1}
\end{equation}
Consider, for $0<\epsilon<1$,
\begin{equation}
v_\epsilon(x):= v(x)+\epsilon d(x),
\label{vepsilon}
\end{equation}
where $d(x):=dist(x, E)$ denotes the distance function from $x$ to $E$.

Then, using (\ref{C5-1}),  there exists some constants $\bar r, c_1, c_2>0$,
depending only on $n$ and $E$,
  such that
\begin{eqnarray*}
\int_{ \partial B_r(0) }(u-v_\epsilon)&=&\int_{ \partial B_r(0)
}(u-v-C_1|x|^2)-\int_{ \partial B_r(0) }(\epsilon d(x)-C_1|x|^2)\\
&\le & -\int_{ \partial B_r(0) }(\epsilon d(x)-C_1|x|^2)\le
-c_1\epsilon r^n +c_2C_1r^{n+1},\quad \forall\ 0<r<\bar r.
\end{eqnarray*}
It follows, for small $\epsilon>0$,  that
\begin{equation}
\lambda(\epsilon):=-\inf_{\Omega\setminus E }(u-v_\epsilon)>0.
\label{negative}
\end{equation}
 Since $v_\epsilon=v$ on $E$, and $u-v>0$ on $\overline \Omega\setminus E$,
we have, for small $\epsilon>0$,
 $$
u-v_\epsilon>0\ \mbox{on}\ \partial \Omega, \qquad
 \liminf_{ x\in \Omega\setminus E, d(x)\to 0}
 (u-v_\epsilon)(x)\ge 0.
 $$
Therefore,
 there
exists $x_\epsilon\in  \Omega\setminus E$ such that
\begin{equation}
u-\tilde v_\epsilon\ge 0\ \mbox{in}\
 \Omega\setminus E, \qquad
(u-\tilde v_\epsilon)(x_\epsilon)=0,
\label{13-1newnew}
\end{equation}
where
 \begin{equation}\tilde v_\epsilon(x):=v_\epsilon(x)-\lambda(\epsilon)
=v(x)+\epsilon d(x)-\lambda(\epsilon).
\label{ZZ11}
\end{equation}
Using the positivity of $u-v$ in $\overline \Omega \setminus E$,
we obtain from the above that
\begin{equation}
\lambda(\epsilon)=-(u-v_\epsilon)(x_\epsilon)=
v(x_\epsilon)-u(x_\epsilon) +\epsilon d(x_\epsilon)\le \epsilon
d(x_\epsilon),\label{12-1new}
\end{equation}
and
$$
\lim_{\epsilon\to 0}
d(x_\epsilon)=0.
$$

We will establish (\ref{13-3new}), with our new $\tilde v_\epsilon$ and
$x_\epsilon$,  for small $\epsilon>0$.
This would lead to a contradiction as shown
in the proof of Theorem \ref{thm1}.

Now we prove (\ref{13-3new}) for small $\epsilon>0$.
By (\ref{12-1new}), we still have (\ref{15-1new})
and (\ref{16-1new}).

By the result in Appendix A,
\begin{eqnarray}
&&\nabla^2 \tilde v_\epsilon(x)-\nabla ^2 v(x)\nonumber\\
&=& \epsilon O(x)^t
\left(diag(0, \cdots, 0, \frac 1{ d(x)}, \cdots, \frac 1{
d(x)})+ O(1)\right) O(x),
\label{CC2}
\end{eqnarray}
where $O(x)$ are orthogonal matrices, and there are $k+1$ zeros
in the diagonal matrix.

Using (\ref{CC2}) instead of (\ref{CC1}),
and using $d(x_\epsilon)\to 0$ instead of $|x_\epsilon|\to 0$, we
establish  (\ref{13-3new})
the same way as in the  proof of Theorem \ref{thm1}.
Theorem \ref{thm11new} for
$k=
n-2$
 is established.

\bigskip

Finally we give the

\noindent{\bf Proof of Theorem \ref{thm11new}.}\
Shrinking $\Omega$ slightly, we may assume that
$u>v$
 on $\overline \Omega\setminus E$,  $u\in C^2(\overline \Omega\setminus E)$,
and $v\in  C^2(\overline \Omega)$.
 We prove the theorem by
contradiction.   Suppose the contrary, then (\ref{zeronewnew})
holds  for some $\bar
x\in E$.  We may assume $\bar x=0\in E$, and
 the tangent space of $E$ at $0$ is spanned
by $e_{n-k+1}, \cdots, e_n$, where $e_1=(1, 0, \cdots, 0)$,
$\cdots, e_n=(0, \cdots, 0, 1)$ are the standard basis of
$\Bbb R^n$.  We write $x=(x', x'')$, where
$x'=(x_1, \cdots, x_{n-k})$ and
$x''=(x_{n-k+1}, \cdots, x_n)$.

For $x$ close to $0$, we have,  for some constant $C$,
\begin{equation}
\frac 34 |x'|-C|x''|^2\le
|x'|-C|x|^2\le
 d(x)\equiv dist(x, E)
\le |x'|+C|x|
\le
\frac 54 |x'|+C|x''|^2.
\label{48-1}
\end{equation}
For example, this follows easily from (\ref{60-1}).

Let $v_\epsilon$ be defined as in (\ref{vepsilon}).  We have

\begin{lem}
For small $\epsilon>0$, (\ref{negative}) holds.
\label{lem-negative}
\end{lem}

\noindent{\bf Proof of Lemma \ref{lem-negative}.}\ We prove it by contradiction.
Suppose not,
then, for some $\bar r, \bar \epsilon, \overline C>0$,
$$
u(x)\ge v_{\bar\epsilon}(x)=
\bar \epsilon d(x)+v(x)\ge
v(0)+\nabla v(0)\cdot x+
\frac {\bar  \epsilon} 2
|x'|-\overline C|x''|^2,\qquad x\in B_{\bar r}\setminus E.
$$

For $0<a<
\min\{ \bar r, \bar \epsilon/4\}$ which will be chosen later, let
$$
h(x):= v(0)+\nabla v(0)\cdot x+
\frac {\bar \epsilon}{ 4a}|x'|^2-(\overline C+1)|x''|^2+
\frac {a^2}2.
$$

On $\partial B_a(0)\setminus E$,
\begin{eqnarray*}
u(x)&\ge &v(0)+\nabla v(0)\cdot x+
  \frac {\bar \epsilon}{2a}|x'|^2
-\overline C|x''|^2 =h(x)+  \frac {\bar \epsilon}{4a}|x'|^2
+|x''|^2-\frac{a^2}2\ge h(x).
\end{eqnarray*}

By Lemma \ref{lemB-2}
in Appendix B and the assumption (\ref{ABC1}), there exists some
positive constant $\bar a>0$ such that if we further
require $0<a<\bar a$, we have
$$
\sum_{i=1}^{k+2}\lambda_i(\nabla ^2 (u-h))
\le \sum_{i=1}^{k+2}\lambda_i(\nabla ^2 u)\le 0.
$$
Indeed, the first inequality follows from
Lemma \ref{lemB-2} by taking
$l=k+2$, $M= \frac {2a}{\bar \epsilon}\nabla^2 u$,
$D=diag(1, \cdots, 1, -\delta_1, \cdots, -\delta_k)$,
and $\delta_1=\cdots=\delta_k= \frac{4a(\overline C+1)}{\bar \epsilon}$;
while the second inequality is (\ref{ABC1}).
Thus, in view of Proposition \ref{prop1},
$$
u-h\ge \inf_{ \partial B_a(0) \setminus E}(u-h)
\ge 0,
$$
and therefore
$$
\liminf_{x\to 0, x\in \Omega\setminus E}u(x)\ge h(0)=v(0)+\frac {a^2}2>v(0),
$$
violating (\ref{zeronewnew}).
Lemma \ref{lem-negative} is established.

\medskip

Given Lemma \ref{lem-negative}, we have, as in the proof of
Theorem \ref{thm11new} for
$k=
n-2$, (\ref{13-1newnew}) for some
$x_\epsilon\in  \Omega\setminus E$, where $\tilde v_\epsilon$ is
given in (\ref{ZZ11}).
The rest of the proof is the same as those in the proof of
Theorem \ref{thm11new} for
$k=
n-2$.
 Theorem \ref{thm11new} is established.

\vskip 5pt
\hfill $\Box$
\vskip 5pt

\section{Proof of Theorem \ref{thmnew1}}

\noindent{\bf Proof of Theorem \ref{thmnew1}.}\

By (\ref{abab0}) and Lemma \ref{lem-4.1-new},
\begin{equation}
\Delta (u-v)\le C-\Delta v\le C_3:=C
-\min_{\overline \Omega}\Delta v<\infty, \qquad
\mbox{in}\ \Omega.
\label{G7}
\end{equation}
We take $u$ the canonical representative in $\Omega$, i.e.
$$
u(x)=\lim_{r\to 0^+}\av_{ B_r(x)  }u,\qquad x\in \Omega,
$$
see Appendix F.

We prove the theorem by contradiction. 
Suppose the contrary, we have, after a possible translation,
$0\in E$, and 
$$
%\begin{equation}
\liminf_{x\in \Omega\setminus E,  x \to 0}(u-v)(x)=0.
$$
%\label{G3}
%\end{equation}
Shrinking $\Omega$ slightly,
we may assume without loss
of generality that $v\in C^2(\overline \Omega)$
and, for some constants $\delta_1, \delta_2>0$,  that
\begin{equation}
(u-v)(x)\ge \delta_1, \qquad x\in \Omega, \  dist(x, \partial \Omega)
\le \delta_2.
\label{G1}
\end{equation}
For $\epsilon\in (0, \delta_1)$,  
let, for $d=diam(\Omega)$, 
$$
w_\epsilon(x):=
\left\{
\begin{array}{ll}
\min\{(u-v)(x)-\epsilon, 0\}& x\in \Omega,\\
0& x\in B_{2d}(0)\setminus \overline \Omega,
\end{array}
\right.
$$
and 
$$
%\begin{equation}
\Gamma_{w_\epsilon} (x):=
\sup\{a+b\cdot x\ |\ a\in \Bbb R, b\in \Bbb R^n,
a+b\cdot z\le w_\epsilon(z)\ \forall\ z\in 
B_{2d}(0)\}
$$
%\end{equation}
be the convex envelope of $w_\epsilon$ on $B_{2d}(0)$.

Since $\Gamma_{w_\epsilon}=0$ outside $
\{x\in \Omega\ |\ dist(x, \partial\Omega)>\delta_2\}$, and
$\min_{ B_{2d}(0) }\Gamma_{w_\epsilon}\le w_\epsilon(0)=
-\epsilon<0$,
the contact set of $w_\epsilon$ and $\Gamma_{w_\epsilon}$ satisfies
\begin{equation}
\{ x\in B_{2d}(0)\ 
|\  w_\epsilon(x)=\Gamma_{w_\epsilon}(x) \}
\subset \{x\in \Omega\ |\ dist(x, \partial\Omega)>\delta_2\}
\subset B_d(0).
\label{G2}
\end{equation}

\begin{lem}  There exist some positive constants $K$
and $\delta$ such that
for any point $\bar x\in \{ x\in B_{2d}(0)\
|\  w_\epsilon(x)=\Gamma_{w_\epsilon}(x) \}$,
there exists $\bar p\in \Bbb R^n$ so that
\begin{equation}
\Gamma_{w_\epsilon}(x) \le 
\Gamma_{w_\epsilon}(\bar x)+\bar p\cdot (x-\bar x)+
K|x-\bar x|^2, \qquad \forall\ |x-\bar x|< \delta.
\label{G6}
\end{equation}
\label{lem-35}
\end{lem}

\noindent{\bf Proof.}\
By (\ref{G1}) and (\ref{G2}),
$
\bar x\in \Omega,$
$dist(\bar x, \partial \Omega)\ge \delta_2,
$
and
$$
u(\bar x)-v(\bar x)-\epsilon =w_\epsilon(\bar x)=
\Gamma_{w_\epsilon}(\bar x)<0,
\quad \mbox{and}\ \ u-v-\epsilon\ge w_\epsilon\ge  \Gamma_{w_\epsilon}\
\mbox{in}\ \Omega.
$$

By (\ref{G7}) and Lemma A in the introduction,
$\Gamma_{w_\epsilon}$ has a unique supporting plane
at $(\bar x, \Gamma_{w_\epsilon}(\bar x))$. So
$\nabla \Gamma_{w_\epsilon}(\bar x)$, as
the slope of the supporting plane, is well 
defined.
Let $\bar y$ be a point near $\bar x$ and let 
$q\in \Bbb R^n$ be the slope of a supporting plane
of $\Gamma_{w_\epsilon}$ at $(\bar y, \Gamma_{w_\epsilon}(\bar y))$, then,
for all $z$ in $\Omega$,
$$
(u-v)(z)\ge (u-v)(\bar x)+\nabla \Gamma_{w_\epsilon}(\bar x)
\cdot (z-\bar x),
$$
and
\begin{eqnarray*}
(u-v)(z)&\ge& \epsilon+\Gamma_{w_\epsilon}(\bar y)
+q\cdot (z-\bar y)\ge  \epsilon+\Gamma_{w_\epsilon}(\bar x)+
\nabla \Gamma_{w_\epsilon}(\bar x)
\cdot (\bar y-\bar x)+q\cdot (z-\bar y)\\
&=&  (u-v)(\bar x)+\nabla \Gamma_{w_\epsilon}(\bar x)
\cdot (\bar y-\bar x)+q\cdot (z-\bar y).
\end{eqnarray*}
Thus, by Lemma \ref{lem-F3} and in view of (\ref{G7}),
$$
|q-\nabla \Gamma_{w_\epsilon}(\bar x)|\le C_4C_3|\bar y|
$$
for some dimensional constant $C_4$.
This implies 
(\ref{G6}).   

\vskip 5pt
\hfill $\Box$
\vskip 5pt

With Lemma \ref{lem-35}, we can apply 
Lemma 3.5 of
\cite{CC} to obtain
 that $\Gamma_{w_\epsilon}\in C^{1,1}_{loc}(B_{2d}(0))$.
Note that in the statement of Lemma 3.5 of
\cite{CC}, the function  is assumed to be continuous,
but the proof there
applies to LSC functions, 
  as in our case. 
By the Alexandrove-Bakelman-Pucci estimate,
\begin{equation}
\epsilon^n=|\inf_ \Omega w_\epsilon|^n
\le \int_{ \{w_\epsilon=\Gamma_{w_\epsilon} \} } \det(\nabla^2 
\Gamma_{w_\epsilon})=
 \int_{ \{w_\epsilon=\Gamma_{w_\epsilon} \}  \cap
(\Omega\setminus E)}  \det(\nabla^2
\Gamma_{w_\epsilon}).
\label{G5}
\end{equation}
For the last equality, we have used the fact that 
$E$ has zero Lebesgue measure.

At a point $x$ in
 $\{w_\epsilon=\Gamma_{w_\epsilon} \} \cap
(\Omega\setminus E)  $,
there exists $a\in \Bbb R$ and $b
$ in $\Bbb R^n$ such that
$$
0>(u-v)(x)-\epsilon=w_\epsilon(x)=\Gamma_{w_\epsilon}(x)\ge -\epsilon,
$$
$$
(u-v)(z)-\epsilon=w_\epsilon(z)\ge
\Gamma_{w_\epsilon}(z)
\ge \Gamma_{w_\epsilon}(x)+\nabla \Gamma_{w_\epsilon}(x)\cdot
(z-x).
$$
If $\nabla \Gamma_{w_\epsilon}$ is differentiable at $x$, then,
by (\ref{abab0}), (\ref{abab5}) and the above,
$$
\Delta  (\Gamma_{w_\epsilon}+v)(x)\le C,
$$
and
$$
 F(x, (\Gamma_{w_\epsilon}+v)(x)+\epsilon,
\nabla (\Gamma_{w_\epsilon}+v)(x),
\nabla ^2 (\Gamma_{w_\epsilon}+v)(x))\ge
F(x, v(x),  \nabla v(x), \nabla^2 v(x)).
$$ 
Clearly,
$$
|\Gamma_{w_\epsilon}(x)|\le  \epsilon,
 \qquad |\nabla \Gamma_{w_\epsilon}(x)|\le \frac
 \epsilon d.
$$
Since $v\in C^2(\overline \Omega)$, 
we have 
$
\Delta  
\Gamma_{w_\epsilon}\le C-\Delta v(x)\le \widehat C$.
By the convexity of $\Gamma_{w_\epsilon}$, $\nabla^2 \Gamma_{w_\epsilon}(x)
\ge 0$, therefore
$$
|\nabla^2 \Gamma_{w_\epsilon}(x)|\le \widehat C
$$
for some positive constant $\widehat C$ independent of $\epsilon$ and $x$.
In the following $\widehat C$ will denote various such constants.

Thus, by (\ref{elliptic}) and the regularity of $F$,
\begin{eqnarray*}
0&\le & F(x, v(x), \nabla v(x), \nabla^2 
(\Gamma_{w_\epsilon}+v)(x))-
F(x, v(x),  \nabla v(x), \nabla^2 v(x))
+\widehat C\epsilon\\
&=&\int_0^1
 \sum_{ij} \frac { \partial F}{ \partial M_{ij} }
(x, v(x),  \nabla v(x), t\nabla^2 \Gamma_{w_\epsilon}(x))
\partial_{ij} \Gamma_{w_\epsilon}(x)dt+\widehat C\epsilon\\
&\le& - \frac 1{ \widehat C}\det(\nabla^2 \Gamma_{w_\epsilon}(x))
+\widehat C\epsilon.
\end{eqnarray*}

By a classical theorem of Rademacher,
$\nabla \Gamma_{w_\epsilon}$ is differentiable almost everywhere.
Therefore, in view of 
(\ref{G5}),
$$
\epsilon^n\le 
\int_{ \{w_\epsilon=\Gamma_{w_\epsilon} \}  \cap
(\Omega\setminus E)}  \det(\nabla^2
\Gamma_{w_\epsilon})
\le \widehat C \epsilon^{n+1}. 
$$
This leads to a contradiction if we choose $\epsilon$ very
small from the beginning.
Theorem \ref{thmnew1} is established.

\vskip 5pt
\hfill $\Box$
\vskip 5pt

\section{  Symmetry of a singular solution }
In this section we give the

\noindent{\bf Proof of Theorem \ref{thm3}.}\
The proof follows closely
arguments in \cite{BN}, with Theorem \ref{thm1} as an
additional ingredient to handle the
possible isolated singularity of $u$ at the origin.

For $0<\lambda<1$, let, with
$x_\lambda:=(2\lambda-x_1, x_2, \cdots, x_n)$,
$$
\Sigma_\lambda:=\{x\ |\ |x|<1, \lambda<x_1<1\},
\qquad
T_\lambda:=\{x\ |\ |x|<1, x_1=\lambda\},
$$
$$
u_\lambda(x):=u(x_\lambda),
\quad\mbox{and}
\quad w_\lambda:= u_\lambda- u.
$$
We will prove that
$$
w_\lambda\ge 0\quad \mbox{in}\
\Sigma_\lambda\setminus\{0_\lambda\}\qquad
\forall\ 0<\lambda<1.
$$

\noindent{\it Step 1.}\
There exists $\frac  12<\lambda_0<1$ such that
$$
w_\lambda\ge 0\quad \mbox{in}\
\Sigma_\lambda\qquad
\forall\ \lambda_0<\lambda<1.
$$
Using   (\ref{cond2}) and (\ref{cond1}),
we have in
$\Sigma_\lambda$, for all $\frac 34 <\lambda< 1$,
\begin{eqnarray}
0&=& F(x, u(x), \nabla u(x), \nabla ^2 u(x))-
F(x_\lambda, u(x_\lambda), \nabla  u(x_\lambda),
\nabla ^2  u(x_\lambda))
\nonumber\\
&\le&  F(x, u(x), \nabla u(x), \nabla ^2 u(x))-F(x, u_\lambda(x),  \nabla u_\lambda(x),
 \nabla ^2  u_\lambda(x)).
\label{order1}
\end{eqnarray}
By (\ref{Q1-1}) and the fact that
$u$ is in $C^2(\overline B_1\setminus\{0\})$,
there exists some positive constant $C_1$ independent of
$\lambda$, and some
functions
$\{a_{ij}(x)\}$, $\{b_i(x)\}$, $c(x)$ satisfying,
with $I$ denoting the $n\times n$ identity matrix,
$$
\frac 1{C_1}I\le (a_{ij}(x))\le C_1I,\ \ \
|b_i(x)|+|c(x)|\le C_1, \qquad \mbox{in}\
\Sigma_\lambda,
$$
such that
\begin{equation}
-a_{ij}(x)\partial_{ij}w_\lambda+b_i(x)\partial_i w_\lambda+c(x)w_\lambda
\ge 0,\qquad \mbox{in}\ \Sigma_\lambda.
\label{32}
\end{equation}
It is clear that
$w_\lambda\ge 0$ on $\partial \Sigma_\lambda$.
Using the maximum principle for
domains of small measure  as in
\cite{BN}, we have
$w_\lambda\ge 0$ in  $\Sigma_\lambda$
if $1-\lambda>0$ is smaller than
some positive constant which depends only on
$C_1$ and $n$.
Step 1 is established.

\medskip

Define
$$
\bar \lambda:=\inf \{ \mu\ |\ 0<\mu<1,
w_\lambda\ge 0\ \mbox{in}\ \Sigma_\lambda\setminus\{0\}\ \forall\
\mu<\lambda<1\}.
$$
Because of Step 1, $0\le \bar \lambda<1$.

\medskip

\noindent{\bf Step 2.}\ $\bar \lambda=0$.

\medskip

There are several cases to consider.

\noindent\underline{Case 1}.\ $\frac 12
<\bar \lambda< 1$.

In this case we argue exactly as in
\cite{BN}. For $\bar \epsilon:=\frac 14 (\frac 12+\bar\lambda)$,
 $\frac 12< \frac 12 \bar\lambda-\bar \epsilon\le
\lambda\le \bar\lambda$,
$w_\lambda$  satisfies a uniformly elliptic inequality
(\ref{32}) in $\Sigma_\lambda$ as before, with ellipticity constants
independent of $\lambda$ --- though
they may depend on $\bar \epsilon$ due to the possible
singularity of $u$ at $\{0\}$.
By continuity,
$w_{\bar \lambda}\ge 0$ in $\Sigma_{\bar \lambda}$.
Since $w_{\bar \lambda}>0$ on $\partial \Sigma_{\bar \lambda}
\cap  B_1$, it follows from the strong maximum principle
that $w_{\bar \lambda}>0$ in $\Sigma_{\bar \lambda}$.
For $0<\delta$ small, let
$D_\delta$ be the set of points in $\Sigma_{\bar \lambda}$
whose distance to its boundary is $\ge \delta$.
Then, in $D_\delta$,
$w_{\bar \lambda}\ge \alpha$ for some positive constant $\alpha$.

For $0<\epsilon<\bar \epsilon$ small,
$w_{  \bar \lambda -\epsilon}\ge \alpha/2$ in $D_\delta$.
Thus
$$
w_{  \bar \lambda -\epsilon}\ge 0\qquad\mbox{on}\
\partial (\Sigma_{  \bar \lambda -\epsilon}\setminus D_\delta).
$$

For $\delta$ and $\epsilon$ small we conclude again,
by using the maximum principle in domains of
small measure as in \cite{BN},
that
$$
w_{  \bar \lambda -\epsilon}\ge 0\ \ \
\mbox{in}\  (\Sigma_{  \bar \lambda -\epsilon}\setminus D_\delta)
\ \mbox{and hence in}\  \Sigma_{\bar \lambda-\epsilon},
$$
contradicting the definition of $\bar \lambda$.

\medskip

\noindent\underline{Case 2}.\ $\bar \lambda=\frac 12$.

Observe first that by property
(\ref{superaffine}),
\begin{equation}
\liminf_{ x\to 0}u(x)\ge a>0.
\label{33}
\end{equation}
By continuity, $w_{\bar \lambda}\ge 0$ in
$\Sigma_{\bar \lambda}$.
Since $u$ is
small near $(1,0,\cdots,0)$, we deduce from
(\ref{33}) that for $0<\rho$ small,
\begin{equation}
w_{\lambda}\ge a/2\ \ \ \mbox{in}\ (B_\rho(0_\lambda)
\cap \Sigma_{\lambda})
\setminus\{0_\lambda\},\qquad
\forall\ \bar\lambda-\rho\le \lambda\le \bar\lambda.
\label{34}
\end{equation}
For $\bar\lambda-\rho\le \lambda\le \bar\lambda$,
$w_\lambda$  satisfies a uniformly elliptic inequality
(\ref{32}) in $\Sigma_\lambda\setminus  B_\rho(0_\lambda)$
 as before, with ellipticity constants
independent of $\lambda$ --- though
they may depend on $\rho$.
By the strong maximum principle,
$w_{  \bar \lambda }>0$ in $\Sigma_{\bar \lambda}
\setminus  B_\rho(0_{\bar \lambda})$.  From now on $\rho$ is fixed.

As in Case 1, for $\delta$ small, consider
$D_\delta=$ set of points in $\Sigma_{\bar \lambda}$
 whose distance to its boundary is
$\ge \delta$.  As before,
$$
w_{  \bar \lambda }>0\ \ \mbox{in}\ D_\delta.
$$
By continuity,
since $D_\delta$ is closed, we have, for small
 $0<\epsilon< \rho$,
$$
w_{  \bar \lambda -\epsilon}>0\ \ \mbox{in}\ D_\delta.
$$
This and (\ref{34}) with $\lambda=\bar \lambda-\epsilon$
imply that $w_{\bar \lambda -\epsilon}\ge 0$ on
$\partial (\Sigma_{\bar \lambda-\epsilon}
\setminus (D_\delta\cup   B_\rho(0_{\bar \lambda -\epsilon})))$.

For $\epsilon, \delta$ small, we have, since
the region have small measure,
$$
\Sigma_{\bar \lambda-\epsilon}
\setminus (D_\delta\cup  B_\rho(0_{\bar \lambda -\epsilon})),
$$
which has small volume, and conclude that
$
w_{  \bar \lambda -\epsilon} \ge 0$ there, and hence
in $\Sigma_{\bar \lambda-\epsilon}\setminus\{0_{\bar \lambda -\epsilon}\}$.
This again contradicts the definition of $\bar \lambda$.

\medskip

\noindent{\it Case 3.}\ $0 <\bar
\lambda<\frac 12$.

Once more,
 $w_{  \bar \lambda }\ge 0$ in
$\Sigma_{\bar \lambda}\setminus\{0_{\bar \lambda}\}$ and
  $w_{  \bar \lambda }>0$ on $\partial \Sigma_{\bar \lambda}
\cap B_1$.
Since $w_{  \bar \lambda }$ is $C^2$ in $\Sigma_{\bar \lambda}\setminus 
\{0_{\bar \lambda}\}$,
and  $w_{\bar \lambda}$ satisfies
a uniformly elliptic
inequality of the form (\ref{32}) on
any compact subset of  $\Sigma_{\bar \lambda}\setminus \{0_{\bar \lambda}\}$,
  by the
strong maximum principle,
 $w_{  \bar \lambda }>0$
in $\Sigma_{\bar \lambda}\setminus \{0_{\bar \lambda}\}$.

Since inequality (\ref{order1}) holds in
$\Sigma_{\bar \lambda}\setminus\{0_{\bar \lambda}\}$, an application of
Theorem
\ref{thm1} yields
$$
\liminf_{x\to 0} w_{  \bar \lambda }(x)>0,
$$
and so, in some ball $B_{2\rho}(0_{\bar \lambda})$, $\rho>0$,
\begin{equation}
w_{\bar \lambda}>\alpha\ \ \mbox{in}\
B_{2\rho}(0_{\bar \lambda})\setminus\{0_{\bar \lambda}\}\ \ \mbox{for some}\
\alpha>0.
\label{99}
\end{equation}

For $0<\lambda\le \bar \lambda$,  $w_ \lambda$ satisfies
a uniformly elliptic
inequality of the form (\ref{32}) in $\Sigma_\lambda\setminus 
B_\rho(0_\lambda)$
with ellipticity constants independent of $\lambda$, though
they may depend on $\rho$.  From now on $\rho$ is fixed.

  As before,
for $0<\delta$ small, let $D_\delta$ be the
set of points in $ \Sigma_{  \bar \lambda }
$ whose distance to its boundary is $\ge \delta$.
We can take $\delta>0$ small so that
the origin $\{0_{\bar \lambda}\}$ is in the
interior of $D_\delta$.   Using (\ref{99}) and the positivity of
$w_{\bar \lambda}$ in $\Sigma_{\bar \lambda}\setminus\{0\}$,
we have
\begin{equation}w_{\bar \lambda}\ge \gamma\ \mbox{in}\
D_\delta \setminus\{0_{\bar \lambda}\}\
\mbox{for some}\ \gamma>0.
\label{98}
\end{equation}
Now consider $\Sigma_{  \bar \lambda -\epsilon}$
for
$0<\bar\lambda-\epsilon$.

By continuity and by (\ref{98}), we have, for   $\epsilon$ small,
$$w_{  \bar \lambda -\epsilon}(x)
>\gamma/2
\ \ \forall\ x\in  D_\delta\setminus\{0_{\bar \lambda-\epsilon}\}.
$$
The above for $x\in (D_\delta\cap B_\rho(0_{\bar \lambda}))
\setminus\{0_{\bar \lambda}\}$ requires some explanation as follows.
We write
$$
w_{  \bar \lambda -\epsilon}(x)=
w_{  \bar \lambda }(\widetilde x)+u(\widetilde x)
-u(x),
$$
where $\widetilde x=(x_1+2\epsilon, x_2, \cdots, x_n)$.
For $\epsilon$ small, $\widetilde x\in B_{2\rho}(0_{\bar \lambda})$.
So, by (\ref{98}),
and by the continuity of $u$ in $\overline D_\delta$, we have,
for small $\epsilon$,
$$
w_{  \bar \lambda -\epsilon}(x)=
w_{  \bar \lambda }(\widetilde x)+u(\widetilde x)
-u(x)\ge \alpha +u(\widetilde x)
-u(x)\ge \alpha/2.
$$

Clearly, $w_{  \bar \lambda -\epsilon}\ge 0$ on
$\partial (\Sigma_{  \bar \lambda -\epsilon}
\setminus D_\delta)$.
As pointed out earlier, $w_{  \bar \lambda -\epsilon}$
 satisfies  a uniformly elliptic inequality like (\ref{32})
in $\Sigma_{  \bar \lambda -\epsilon}
\setminus D_\delta$
with ellipticity constants independent of $\epsilon$ and $\delta$.
Applying the maximum principle
in 
the region,
 which
has small measure,
 we conclude that
$w_{  \bar \lambda -\epsilon}\ge 0$ there, and hence in
$\Sigma_{  \bar \lambda -\epsilon}\setminus\{0_{\bar\lambda -\epsilon}\}$.  Contradiction.

Step 2 is established.

Since $w_\lambda=0$ on $T_\lambda$, it follows from Step 2 that
$$
2
\frac { \partial u}{ \partial x_1}=
 \frac { \partial w_\lambda}{ \partial x_1}\le 0,\qquad
\mbox{on}\ T_\lambda,\ \forall\ 0<\lambda<1,
$$
namely,
$$
\frac { \partial u}{ \partial x_1}(x)\le 0\qquad
\forall\ x\in B_1, \ 0<x_1<1.
$$

For $0<\lambda<1$, we
 already know that $w_\lambda\ge 0$ in $\Sigma_\lambda
\setminus\{0_\lambda\}$.
We also know that $w_\lambda$ is not identically zero since
it is positive on $\partial \Sigma_\lambda\cap  B_1$.
Because of (\ref{order1}), we can apply the Hopf Lemma
to  $w_\lambda$ to obtain
$$
 \frac { \partial w_\lambda}{ \partial x_1}<0,\qquad
\mbox{on}\ T_\lambda,\ \forall\ 0<\lambda<1,
$$
namely,
 (\ref{mono}) holds.
Theorem \ref{thm3} is established.

\section{Proof of Theorem \ref{thm2-1new}}

In  this section we give the

\medskip

 \noindent{\bf Proof of Theorem
\ref{thm2-1new}.}\ It is similar to that of
 Theorem \ref{thm11new}.
Shrinking $\Omega$ slightly, we may assume that
$u>v$
 on $\overline \Omega\setminus \{0\}$,
 $u\in C^2(\overline \Omega\setminus \{0\})$,
and $v\in  C^2(\overline \Omega)$.
 We prove the theorem by
contradiction.   Suppose the contrary, then 
(\ref{zeronew})
holds.

Fix some $0<\alpha<1$.  Let
$$
v_\epsilon(x)=v(x)+\epsilon |x|^{1+\alpha}.
$$
Similar to Lemma \ref{lem-negative} we have
\begin{lem}
For  $\epsilon>0$, 
$$
%\begin{equation}
\lambda(\epsilon):=-\inf_{\Omega\setminus \{0\}  }(u-v_\epsilon)>0.
$$
%\label{41new}
%\end{equation}
\label{lem-D1}
\end{lem}

\noindent{\bf Proof.}\  We prove it by contradiction.
Suppose not,
then, for some $ \epsilon>0$ and $C>0$,
$$
u(x)\ge v_\epsilon(x)
\ge \epsilon |x|^{1+\alpha}+v(0)+\nabla v(0)\cdot x-C|x|^2,
\quad x\in \Omega\setminus\{0\}.
$$
For small $r>0$,
$$
\min_{ x\in \partial B_r}[u(x)-\nabla v(0)\cdot x]
\ge v(0)+\epsilon r^{1+\alpha}-Cr^2\ge v(0)+
\frac \epsilon 2 r^{1+\alpha}.
$$

On the other hand,
in view of 
(\ref{zeronew}),
$$
\inf_{ x\in B_r\setminus\{0\} }
[u(x)-\nabla v(0)\cdot x]
\le \liminf_{ x\to 0} [u(x)-\nabla v(0)\cdot x]
=v(0)<  v(0)+
\frac \epsilon 2 r^{1+\alpha}.
$$

The above contradicts to condition  (\ref{superaffine}).

\vskip 5pt
\hfill $\Box$
\vskip 5pt

By Lemma \ref{lem-D1},
 there
exists $x_\epsilon\in  \Omega\setminus \{0\}$  such that
$$
u-\tilde v_\epsilon\ge 0\ \mbox{in}\
 \Omega\setminus \{0\}, \qquad
(u-\tilde v_\epsilon)(x_\epsilon)=0,
$$
where
$$
\tilde v_\epsilon(x):=v_\epsilon(x)-\lambda(\epsilon)
=v(x)+\epsilon  |x|^{1+\alpha}-\lambda(\epsilon).
$$

Using the positivity of $u-v$ in $\overline \Omega \setminus  \{0\}$,
we obtain
\begin{equation}
0<\lambda(\epsilon)=-(u-v_\epsilon)(x_\epsilon)=
v(x_\epsilon)-u(x_\epsilon) +\epsilon |x_\epsilon|^{1+\alpha}
\le \epsilon
 |x_\epsilon|^{1+\alpha},\label{12-1newnewnew}
\end{equation}
and
$$
\lim_{\epsilon\to 0}
 |x_\epsilon|=0.
$$

We will show that
\begin{equation}
\widetilde A^{ \tilde v_\epsilon}(x_\epsilon)\in
{\cal S}^{n\times n} \setminus
U. \label{13-4new}
\end{equation}
A calculation gives
\begin{eqnarray}
&&\nabla^2 \tilde v_\epsilon(x)-\nabla ^2 v(x)\nonumber\\
&=&
(1+\alpha) \epsilon O(x)^t
\left(diag( \alpha |x|^{\alpha-1},
|x|^{\alpha-1}, \cdots, |x|^{\alpha-1})+ O(1)\right) O(x),
\nonumber
%\label{CC2newnew}
\end{eqnarray}
where $O(x)$ are orthogonal matrices.

At $x_\epsilon$, using (\ref{12-1newnewnew}),
$$
|L(x_\epsilon, v, \nabla v)-L(x_\epsilon, \tilde v_\epsilon, \nabla
\tilde v_\epsilon)| \le  C(|\nabla v-\nabla\tilde
v_\epsilon|+|v-\tilde v_\epsilon|) \le C\epsilon
|x_\epsilon|^\alpha.
$$
It follows, for $\epsilon>0$ small,  that
$$
\widetilde A^v(x_\epsilon) \ge \widetilde A^{\tilde
v_\epsilon}(x_\epsilon)+\alpha(1+\alpha)\epsilon
|x_\epsilon|^{\alpha-1}I-C\epsilon |x_\epsilon|^\alpha I> \widetilde A^{\tilde
v_\epsilon}(x_\epsilon).
$$
Since $\widetilde A^v(x_\epsilon)$  is in ${\cal S}^{n\times n}
\setminus U$, so is $\widetilde A^{\tilde
v_\epsilon}(x_\epsilon)$ ----  due to
 (\ref{YY1}).

Since $x_\epsilon$ is an interior local minimum
point of $u-\tilde v_\epsilon$, we have
$$
u(x_\epsilon)=\tilde v_\epsilon(x_\epsilon),\ \nabla
u(x_\epsilon)=\tilde v_\epsilon(x_\epsilon),\ \nabla
^2u(x_\epsilon)\ge \nabla ^2 v_\epsilon(x_\epsilon).
$$
Thus
$$
\widetilde A^{\tilde v_\epsilon}(x_\epsilon)\ge \widetilde
A^u(x_\epsilon). 
$$
Since  $\widetilde A^u(x_\epsilon)\in\overline U$, we
know from the above, in view of  (\ref{YY1}), that
 $\widetilde A^{\tilde
v_\epsilon}(x_\epsilon)\in \overline U$.  This violates
(\ref{13-4new}).  Theorem \ref{thm2-1new} is established.

\section{Proof of Theorem \ref{thm2-1}}

In  this section we give the

\medskip

 \noindent{\bf Proof of Theorem
\ref{thm2-1}.}\ It uses arguments in the proofs of
 Theorem \ref{thm11new} and Theorem \ref{thm2-1new}.

We start with the first paragraph of the
proof of  Theorem \ref{thm11new}.
Fix some constant $0<\alpha<1$.  For  $x$ close to $0$, we have,  for some constant $C$,
$$
\frac 34 |x'|^{1+\alpha}-C|x''|^2\le d(x)^{1+\alpha}
\equiv dist(x, E)^{1+\alpha}\le
\frac 54 |x'|^{1+\alpha}+C|x''|^2.
$$
We only prove the left inequality.  From (\ref{48-1}) we have
$$
|x'|\le d(x)+C|x''|^2,
$$
so
\begin{eqnarray*}
|x'|^{1+\alpha}&\le &
d(x)^{1+\alpha} \left( 1+C(\frac{|x''|^2} {d(x)})\right)^{1+\alpha}\\
&\le &
d(x)^{1+\alpha}  \left( 1+ K \frac{|x''|^2} {d(x)}
+K(\frac{|x''|^2} {d(x)})^{1+\alpha}\right)\\
&\le &d(x)^{1+\alpha} +Kd(x)^\alpha|x''|^2
+K |x''|^{ 2(1+\alpha) }
\end{eqnarray*}
which yields the first inequality.

Let
$$
v_\epsilon(x)=v(x)+\epsilon d(x)^{1+\alpha}.
$$

\noindent{\bf Claim.}\
Lemma \ref{lem-negative} holds with the above $v_\epsilon$.

\medskip

To prove the claim, we only need to follow the proof of
Lemma \ref{lem-negative} with modification as below:

1. In the first displayed formula, change ``$\bar \epsilon d(x)$''
to ``$\bar \epsilon d(x)^{1+\alpha}$'', and change
``$\frac {\bar  \epsilon} 2
|x'|$'' to
``$\frac {\bar  \epsilon} 2
|x'|^{1+\alpha}$''.

2.  Change ``$0<a<
\min\{ \bar r, \bar \epsilon/4\}$''
to  ``$0<a<
\min\{ \bar r, ( \bar \epsilon/4)^{1-\alpha})\}$''.

3. In the next two displayed formula,
change the two ``$\frac {\bar \epsilon}{ 4a}$''
to  ``$\frac {\bar \epsilon}{ 4 a^{1-\alpha}}$'',
and one  ``$\frac {\bar \epsilon}{2a}$''
to  ``$\frac {\bar \epsilon}{2 a^{1-\alpha}}$''.

4. In the 4th and 5th lines below, change
  ``$M= \frac {2a}{\bar \epsilon}\nabla^2 u$''
to  ``$M= \frac {2a^{1-\alpha}}{\bar \epsilon}\nabla^2 u$'',
and change
``$\delta_1=\cdots=\delta_k= \frac{4a(\overline C+1)}{\bar \epsilon}$''
to ``$\delta_1=\cdots=\delta_k=
\frac{4a^{1-\alpha}(\overline C+1)}{\bar \epsilon}$''.

\medskip

Given the claim above,  we have, as in the proof of
Theorem \ref{thm11new} for
$k=
n-2$, (\ref{13-1newnew}) for some
$x_\epsilon\in  \Omega\setminus E$, where
$$\tilde v_\epsilon(x):=v_\epsilon(x)-\lambda(\epsilon)
=v(x)+\epsilon d(x)^{1+\alpha}-\lambda(\epsilon).
$$

Using the positivity of $u-v$ in $\overline \Omega \setminus E$,
we obtain
\begin{equation}
0<\lambda(\epsilon)=-(u-v_\epsilon)(x_\epsilon)=
v(x_\epsilon)-u(x_\epsilon) +\epsilon d(x_\epsilon)^{1+\alpha}
\le \epsilon
d(x_\epsilon)^{1+\alpha},\label{12-1newnew}
\end{equation}
and
$$
\lim_{\epsilon\to 0}
d(x_\epsilon)=0.
$$

We will show that
$$
%\begin{equation}
\widetilde A^{ \tilde v_\epsilon}(x_\epsilon)\in
{\cal S}^{n\times n} \setminus
U.
$$
% \label{13-4}
%\end{equation}

Since $x_\epsilon$ is an interior local minimum
point of $u-\tilde v_\epsilon$, we have
$$
u(x_\epsilon)=\tilde v_\epsilon(x_\epsilon),\ \nabla
u(x_\epsilon)=\tilde v_\epsilon(x_\epsilon),\ \nabla
^2u(x_\epsilon)\ge \nabla ^2 v_\epsilon(x_\epsilon).
$$
Thus
$$
%\begin{equation}
\widetilde A^{\tilde v_\epsilon}(x_\epsilon)\ge \widetilde
A^u(x_\epsilon).
$$
% \label{17-2}
%\end{equation}

By the result in Appendix A,
\begin{eqnarray}
&&\nabla^2 \tilde v_\epsilon(x)-\nabla ^2 v(x)\nonumber\\
&=&
(1+\alpha) \epsilon O(x)^t
\left(diag(0, \cdots, 0, \alpha d(x)^{\alpha-1},
d(x)^{\alpha-1}, \cdots, d(x)^{\alpha-1})+ O(1)\right) O(x),
\label{CC2new}
\end{eqnarray}
where $O(x)$ are orthogonal matrices, and there are $k$ zeros
in the diagonal matrix.

At $x_\epsilon$, using (\ref{12-1newnew}),
$$
|L(x_\epsilon, v, \nabla v)-L(x_\epsilon, \tilde v_\epsilon, \nabla
\tilde v_\epsilon)| \le  C_1(|\nabla v-\nabla\tilde
v_\epsilon|+|v-\tilde v_\epsilon|) \le C_1\epsilon
d(x_\epsilon)^\alpha,
$$
where $C_1>0$ is some constant independent of $\epsilon$.
In the following we use $\{C_i\}$ to denote various 
positive constants which are  independent of $\epsilon$.

Since $d(x_\epsilon)\to 0$ and $\alpha<1$,
we have
\begin{equation}
\widetilde A^v(x_\epsilon) -\widetilde A^u(x_\epsilon)
\ge \widetilde A^v(x_\epsilon) - \widetilde A^{\tilde
v_\epsilon}(x_\epsilon)\ge \hat N,
\label{WWW1}
\end{equation}
where
$$
\hat N:=\alpha
 \epsilon O(x_\epsilon)^t
\left(diag(-C_2, \cdots, -C_2,  d(x_\epsilon)^{\alpha-1},
d(x_\epsilon)^{\alpha-1}, \cdots, d(x_\epsilon)^{\alpha-1})\right) 
O(x_\epsilon).
$$
Since $v$ is $C^2$ in $\Omega$, 
$$
%\begin{equation}
M:=\widetilde A^u(x_\epsilon)\le
\widetilde A^v(x_\epsilon) -\hat N\le \widetilde A^v(x_\epsilon)+ 
\alpha C_2\epsilon I\le
 C_3 I
$$
%\label{WWW2}
%\end{equation}

Let $\bar\delta=\bar \delta(C_3)$ and $\bar \epsilon=\bar \epsilon(C_3)$
 be the positive numbers in Property $P_k$,
and let
$$
N:= \frac{ \alpha
 \epsilon \bar\delta}{ C_2}  O(x_\epsilon)^t
\left(diag(-\bar\delta, \cdots, -\bar\delta,    
 1,  \cdots,    
 1)\right)
O(x_\epsilon).
$$

Since $M=\widetilde A^u(x_\epsilon)\in\overline U$, we know from
 Property $P_k$ that 
$M+N\in U$ for small $\epsilon>0$.

For small $\epsilon>0$,
$ d(x_\epsilon)^{\alpha-1}>
C_2/\bar \delta$,
so $
\hat N\ge N.
$
 
By (\ref{WWW1}) and the above, $\widetilde A^v(x_\epsilon)
\ge M+\hat N\ge M+N$.
Since $M+N\in U$ for small $\epsilon$, we have,
 and in view of  (\ref{YY1}),  
$\widetilde A^v(x_\epsilon)\in U$ for small $\epsilon$.
This violates (\ref{2-2}).
Theorem \ref{thm2-1} is established.

\vskip 5pt
\hfill $\Box$
\vskip 5pt

\noindent{ \bf  About Remark \ref{rem-13}.}\

For (ii):  Let
${\cal N}\subset {\cal S}^{n\times n}$ denote the set
of matrices whose eigenvalues consist
of $k$ zeros and $n-k$ $1^{'}$s.
We know that ${\cal N}\subset U$.  Since 
${\cal N}$ is compact, there exists some positive number
$\bar \delta$ such that $n\bar\delta-$neighborhood of
${\cal N}$ belongs to $U$.  Consequently,
$N\in U$ if its eigenvalues consist of
$k$ $-\bar\delta^{'}$s and $n-k$ $1^{'}$s.
For $M\in \overline U$, we know that $2M\in \overline U$,
since $U$ is invariant under multiplication by positive
constants.
Similarly,  $2\epsilon N\in U$
for any $0<\epsilon<1$.  By the convexity of $U$,
$M+\epsilon N=\frac 12(2M+2\epsilon N)\in U$.

\medskip

For (iii):  It is easy to see that such $U_l$ satisfies
(ii).

\medskip

For (i): 
For any $M\in \partial U$, $M\le CI$, there exists
some $a=a(M)>0$ and $\beta=\beta(M)>0$ such that
if $\hat M\in \partial U$, $|\hat M-M|<\beta$, 
then
$$
\hat M+a O^t diag(-\delta, \cdots, -\delta,
1, \cdots, 1) O\in U,\qquad \forall \ 0<\delta\le \beta,
$$
where there are $k$ $-\delta^{'}$s and
$n-k$ $1^{'}$s, and $O\in O(n+1)$ is an orthogonal matrix.
By the convexity of $U$, we have
$$
\hat M+
\epsilon  O^t diag(-\delta, \cdots, -\delta,
1, \cdots, 1) O\in U,\quad \forall\ 
0<\delta\le \beta, 0<\epsilon\le a.
$$
By the compactness of $\partial U\cap \{M\le CI\}$, there
exist $\bar\delta>0$ and $\bar \epsilon>0$ such that
$$
M+
\epsilon  O^t diag(-\delta, \cdots, -\delta,
1, \cdots, 1) O\in U,\quad \forall\ 
0<\delta\le \bar\delta, 0<\epsilon\le \bar\epsilon.
$$

\section{Appendix A}

For $n\ge 1$ and $0\le k\le n-1$, let
$E$ be a smooth compact  $k-$dimensional
submanifold in $ \Bbb R^n$, and let
$$
d(x):=dist(x, E):=\inf\{ |x-y|\ |\ y\in E\}, \qquad x\in \Bbb R^n,
$$
denote the distance function from $x$ to $E$.
For $x\in \Bbb R^n\setminus E$ close to $E$,
we use $\nabla^2 d(x)$ to denote the Hessian of
$d(x)$.

\begin{lem}
Let $E\subset \Bbb R^n$ be as above, and let $G$ be a smooth
one variable function, $\widehat G(x)= G(d(x))$.
 Then
 there exists some  matrix valued function
 $O(x)\in O(n)$
such that
\begin{eqnarray}
O(x)^t (\nabla  ^2 \widehat G(x)) O(x)
&=& diag\left(0, \cdots, 0, G''(d(x)),  d(x)^{-1}G'(d(x)),
\cdots, d(x)^{-1}G'(d(x))\right)\nonumber\\
&&
+O\left(d|G''(d)|+|G'(d)|\right),
\quad\mbox{as}\ d(x)\to 0, \label{zzzz} \end{eqnarray}
with  $k$ zeros in the
 diagonal matrix.
\label{lemapp1}
\end{lem}

\noindent{\bf Proof.}\
For $x$ close to $E$, there exists a unique $\bar x=
\bar x(x)\in E$ such that
$d(x)=|x-\bar x|$.
Without loss of generality, we may assume that $\bar x=0$.
Then $d(x)=|x|$.
We may also  assume that
$T_0E$, the tangent space of $E$ at $0$,
is spanned by $e_{n-k+1}, \cdots, e_n$, where $e_1=(1, 0, \cdots, 0)$,
$\cdots, e_n=(0, \cdots, 0, 1)$ are the standard basis of
$\Bbb R^n$.
Fixing a smooth local orthonormal frame
 $\{\epsilon_1(p), \cdots, \epsilon_{n-k}(p)\}$
of the orthogonal complement of
$T_pE$, the tangent space of $E$ at $p\in E$,
there is a unique representation of $x\in \Bbb R^n$ near $0$
given by
$$
x=p(x)+\varphi_1(x)\epsilon_1(p(x))+\cdots+
\varphi_{n-k}(x)\epsilon_{n-k}(p(x))
$$
where $\varphi_\alpha$ and $p$ are smooth functions
near $0$.

We may also assume  that
$\epsilon_\alpha(0)=e_\alpha$ for $1\le \alpha\le n-k$.
Clearly
\begin{equation}
\varphi_\alpha(0)=0,
\qquad \nabla \varphi_\alpha(0)=e_\alpha, \qquad 1\le \alpha\le n-k,
\label{app3}
\end{equation}

Near $0$,
\begin{equation}
d(x)=\sqrt{ \sum_{\alpha=1}^{n-k} \varphi_\alpha(x)^2 }.
\label{60-1}
\end{equation}
By chain rule,
$$
\partial_{ij}\widehat G(x)=
G''(d(x))\partial_i d(x) \partial_j d(x)+
G'(d(x)) \partial_{ij} d(x).
$$
A calculation gives
$$
\partial_i d=d^{-1} \sum_\alpha \varphi_\alpha (\partial_i
\varphi_\alpha),
$$
$$
\partial_{ij}d =-d^{-3} \sum_{\alpha, \beta}  \varphi_\alpha
 \varphi_\beta  (\partial_i
\varphi_\alpha)  (\partial_j\varphi_\beta)
+ d^{-1}  \sum_{\alpha}  (\partial_i
\varphi_\alpha)(\partial_j  \varphi_\alpha)
+d^{-1}  \sum_{\alpha}  \varphi_\alpha
(\partial_{ij}  \varphi_\alpha).
$$

Using (\ref{app3}) and
the fact $|\varphi_\alpha(x)|\le d(x)$,
we have,    as $|x|\to 0$, that
$$
\partial_i d(x) \partial_j d(x) = O(d(x)),\qquad \mbox{if}\ \max\{i, j\}>n-k,
$$
$$
\partial_i d(x)
\partial_j d(x) = d^{-2} \varphi_i\varphi_j
+ O(d(x)), \qquad \qquad 1\le i, j\le n-k,
$$
$$
 d(x) \partial_{ij}d(x)=
O(d(x)),\qquad \mbox{if}\ \max\{i, j\}>n-k,
$$
\begin{eqnarray}
 d(x) \partial_{ij}d(x)&=&
 \sum_{\alpha}(\partial_i \varphi_\alpha)
(\partial_j\varphi_\alpha)-
d^{-2} \sum_{\alpha, \beta}(\varphi_\beta \partial_j \varphi_\beta)
(\varphi_\alpha \partial_i\varphi_\alpha)+ \sum_{\alpha}
\varphi_\alpha (\partial_{ij}\varphi_\alpha)\nonumber\\
&=& \delta_{ij}-d^{-2}  \varphi_i
\varphi_j+
O(d(x)),
\qquad \qquad 1\le i, j\le n-k.
\nonumber
\end{eqnarray}
Since the $(n-k)\times (n-k)$ matrix
$(d^{-2}  \varphi_i
\varphi_j)$ has eigenvalues $0$ and $1$,
with $0$ having the multiplicity $n-k-1$,
(\ref{zzzz}) follows immediately.

\section{Appendix B}

\begin{lem} Let $M\in {\cal S}^{n\times n}$, $1\le k\le n$,
then
$$
\sum_{i=1}^k \lambda_i(M)=
\min\{ \sum_{i=1}^k \epsilon_i^t M\epsilon_i\ |\
\epsilon_i\in \Bbb R^n, \epsilon_i^t \epsilon_j =\delta_{ij}\}.
$$
\label{lemB-1}
\end{lem}

\noindent{\bf Proof.}\ There exists
orthonormal basis $\{E_i\}_{i=1}^n$ such that
$$
M E_i =\lambda_i(M) E_i.
$$
Clearly
\begin{equation}
 \sum_{i=1}^k E_i^t M E_i=
 \sum_{i=1}^k \lambda_i(M)  E_i^t  E_i= \sum_{i=1}^k \lambda_i(M).
\label{B-1}
\end{equation}
Let $\{\epsilon\}_{i=1}^k$ satisfy $ \epsilon_i^t \epsilon_j =\delta_{ij}$.
There exists $O\in O(n)$ such that
$$
\epsilon_i=\sum_{j=1}^n O_{ij} E_j, \qquad 1\le i\le k.
$$
Thus
\begin{eqnarray*}
&&\sum_{i=1}^k \epsilon_i^t M\epsilon_i\\
&=&
\sum_{i=1}^k \sum_{j=1}^n  \sum_{l=1}^n
O_{ij}E_j^t O_{il}ME_l=\sum_{i=1}^k \sum_{j=1}^n  \sum_{l=1}^n
O_{ij}O_{il} \lambda_l(M)E_j^t E_l
=
\sum_{i=1}^k \sum_{j=1}^n
O_{ij}^2 \lambda_j(M)\\
&\ge&
\sum_{j=1}^k  \lambda_j(M)
(\sum_{i=1}^k O_{ij}^2)+
\lambda_k(M) (\sum_{j=k+1}^n \sum_{i=1}^k  O_{ij}^2)
\end{eqnarray*}
Since $O\in O(n)$,
$$
 \sum_{j=1}^n  \sum_{i=1}^k O_{ij}^2= k,\qquad \sum_{i=1}^k  O_{ij}^2
\le 1.
$$
It follows that
\begin{eqnarray*}
&&\sum_{i=1}^k \epsilon_i^t M\epsilon_i\\
&\ge&
\sum_{j=1}^k \lambda_j(M)+
\sum_{j=1}^k  \lambda_j(M)
(\sum_{i=1}^k O_{ij}^2-1)
+
\lambda_k(M) (k- \sum_{j=1}^k \sum_{i=1}^k  O_{ij}^2)\\
&\ge & \sum_{j=1}^k \lambda_j(M)+
\sum_{j=1}^k  \lambda_k(M)
(\sum_{i=1}^k O_{ij}^2-1)
+
\lambda_k(M) (k- \sum_{j=1}^k \sum_{i=1}^k  O_{ij}^2)
=\sum_{j=1}^k \lambda_j(M).
\end{eqnarray*}
Lemma \ref{lemB-1} follows from
this and
(\ref{B-1}).

\bigskip

\begin{lem}  Let $n\ge 2$, $1\le k<l\le n$,
 $0<\delta_i\le (l-k)/k$, $1\le i\le k$,
and  $M\in {\cal S}^{n\times n}$.  Then
$$
\lambda_1(M-D)+ \cdots+ \lambda_{l}(M-D)
\le \lambda_1(M)+ \cdots+ \lambda_{l}(M),
$$
where $D=diag(1, \cdots, 1, -\delta_1, \cdots, -\delta_k)$.
\label{lemB-2}
\end{lem}

\noindent{\bf Proof.}\  There exists orthonormal basis
$\{E_i\}$ of $\Bbb R^n$ such that
$$
ME_i=\lambda_i(M)E_i, \qquad 1\le i\le n.
$$
By Lemma \ref{lemB-1},
\begin{equation}
\sum_{i=1}^l \lambda_i(M-D)
\le
\sum_{i=1}^l E_i^t (M-D) E_i=
\sum_{i=1}^l  \lambda_i(M)-
\sum_{i=1}^l E_i^t D E_i.
\label{B-3}
\end{equation}
Let $e_1=(1, 0, \cdots)$,
$\cdots$, $e_n=(0, \cdots, 0, 1)$ be the standard basis.
We know that
$$
De_m=e_m, \ 1\le m\le n-k; \qquad
De_m= -\delta_{m-(n-k)}e_m, \ n-k+1\le m\le n.
$$
There exists $O\in O(n)$ such that
$$
E_i=\sum_{j=1}^n O_{ij}e_j, \qquad 1\le i\le n.
$$
It follows that
\begin{eqnarray*}
\sum_{i=1}^l E_i^t D E_i&=& \sum_{i=1}^l\sum_{j=1}^n
\sum_{m=1}^n O_{ij} e_j^t O_{im} D e_m
\\
&=&  \sum_{i=1}^l\sum_{j=1}^n
\sum_{m=1}^{n-k}  O_{ij}  O_{im} \delta_{jm}-
 \sum_{i=1}^l\sum_{j=1}^n
\sum_{m=n-k+1}^n   O_{ij}  O_{im} \delta_{ m-(n-k)}  \delta_{jm}\\
&=&  \sum_{i=1}^l\sum_{m=1}^{n-k} O_{im}^2
-  \sum_{i=1}^l \sum_{m=n-k+1}^n  O_{im}^2
 \delta_{ m-(n-k)}
\end{eqnarray*}
Since $O\in O(n)$, we have
$$
\sum_{m=1}^n O_{im}^2=1, \qquad  \sum_{i=1}^l O_{im}^2\le 1.
$$
Thus
\begin{eqnarray*}
\sum_{i=1}^l E_i^t D E_i&=&  \sum_{i=1}^l\sum_{m=1}^n  O_{im}^2
-  \sum_{i=1}^l\sum_{m=n-k+1}^n O_{im}^2
-  \sum_{i=1}^l\sum_{m=n-k+1}^n O_{im}^2 \delta_{   m-(n-k) }\\
&\ge & l-k- (\frac {l-k}k)k=0.
\end{eqnarray*}
Lemma \ref{lemB-2} follows from this and (\ref{B-3}).

\section{Appendix C. Proof of Proposition \ref{prop1}}

For the sake of exposition, we  first prove the result for $k=-1$,
then for $k=0$, and finally for $k\ge 1$.

\medskip

\noindent\underline{The proof for $k=-1$}.\
For any $\epsilon>0$, let $u_\epsilon(x):=u(x)-\epsilon |x|^2$.
Since
$\lambda_1(\nabla ^2 u_\epsilon)=\lambda_1(\nabla ^2 u)-2\epsilon
\le -2\epsilon<0$, we easily see that
the minimum of $u_\epsilon$ can only be attained on
$\partial \Omega$.  Indeed, if  $x_\epsilon$ is a minimum point
of $u_\epsilon$ in $\Omega$, then
$\nabla^2 u_\epsilon(x_\epsilon)$ is a semi-positive definite matrix,
violating $\lambda_1(\nabla ^2 u_\epsilon)<0$.
Sending $\epsilon$ to $0$ yields
$\min_{\overline \Omega}u =
\min_{\partial  \Omega}u$.

\medskip

\noindent\underline{The proof for $k=0$}.\
For simplicity we assume that $E$ contains only
one point, say $\{0\}$,  in $\overline \Omega$, since the proof
requires only some modification
 when $E$ consists of finitely many points; see the proof for $k\ge 1$.
Let
$$
\hat u:= u- \inf_{\Omega\setminus\{0\}}u.
$$
Then
\begin{equation}
\inf_{\Omega\setminus\{0\}}\hat u=0.
\label{ABC6}
\end{equation}
We need to prove that
 $$
\inf_{\partial \Omega\setminus\{0\}}\hat u=0.
$$
We prove it by contradiction.  Suppose the contrary, then
 $$
\inf_{\partial \Omega\setminus\{0\}}\hat u=\delta>0.
$$
Let $\bar r=6 diam(\Omega)$.
For any $\epsilon>0$ small, let
$$
h_\epsilon(r):=
1-\frac { \log(\frac r{\bar r})  }
        {   \log(\frac \epsilon{\bar r})  },
\qquad
\epsilon\le r\le \bar r; \qquad\hat h_\epsilon(x):=h_\epsilon(|x|).
$$
Clearly,
$$
h_\epsilon(\epsilon)=0,\qquad h_\epsilon(r)\le 1\ \forall\ \epsilon<r<\bar r,
$$
$$
h_\epsilon'(r)=
-\frac 1{   r  \log(\frac \epsilon{\bar r})   }>0,\qquad
\epsilon<r<  \bar r,
$$
$$
h_\epsilon''(r)=
\frac 1{  r^2  \log(\frac \epsilon{\bar r})   }<0,\qquad
\epsilon<r<  \bar r,
$$
$$
\frac { h_\epsilon'(r) }r =
-\frac 1{  r^2  \log(\frac \epsilon{\bar r})   }>0,\qquad
\epsilon<r<  \bar r.
$$
For $x\in \Omega\setminus B_\epsilon$,
 the $n$ eigenvalues of
$\nabla^2 \hat  h_\epsilon$ are, with $r=|x|$,
$$
h_\epsilon''(r)< \frac { h_\epsilon'(r) }r
=\cdots= \frac { h_\epsilon'(r) }r,
$$
i.e.
$$
\lambda_1(\nabla^2 \hat  h_\epsilon(x))=
h_\epsilon''(r)<0, \lambda_2(\nabla^2 \hat h_\epsilon(x))= \cdots=
 \lambda_n(\nabla^2 \hat h_\epsilon(x))=\frac { h_\epsilon'(r) }r =
-h_\epsilon''(r)>0.
$$
Since
$$
\lambda_1(\nabla^2
(\delta \hat  h_\epsilon))+\lambda_2(\nabla^2
(\delta \hat  h_\epsilon))= 0\qquad
\mbox{in}\
\Omega\setminus  \overline  B_\epsilon,
$$
and
$$
\lambda_1(\nabla^2 \hat u)+\lambda_2(\nabla ^2 \hat  u)\le 0\qquad
\mbox{in}\
\Omega\setminus  \overline  B_\epsilon,
$$
we have,
$$
\lambda_1(\nabla ^2(\hat u-\delta  \hat  h_\epsilon))
=\lambda_1( \nabla ^2 \hat u - \nabla ^2 (\delta  \hat  h_\epsilon)) \le 0,
\qquad
\mbox{in}\
\Omega\setminus  \overline  B_\epsilon.
$$
Indeed if the above does not hold, then
$\nabla ^2 \hat u - \nabla ^2 (\delta  \hat  h_\epsilon)$ is a
positive definite matrix, and therefore
$\lambda_i(\nabla ^2 \hat u)>\lambda_i(\nabla ^2 (\delta  \hat  h_\epsilon))$
for every $i$.
In particular,
$\lambda_1(\nabla^2 \hat u)+\lambda_2(\nabla ^2 \hat  u)
>\lambda_1(\nabla^2
(\delta \hat  h_\epsilon))+\lambda_2(\nabla^2
(\delta \hat  h_\epsilon))$, which leads to contradiction.

Thus, by Proposition \ref{prop1} for $k=-1$,
$$
\hat u-\delta  \hat  h_\epsilon\ge \inf_{\partial (\Omega\setminus B_\epsilon) }
(\hat u-\delta  \hat  h_\epsilon)\ge 0,
\qquad \qquad
\mbox{in}\
\Omega\setminus  \overline  B_\epsilon.
$$
Now, for any $x\in \Omega\setminus \{0\}$,
we have, from the above, that
$$
\hat u(x)\ge \delta  \hat  h_\epsilon(x),
\qquad \forall\ 0<\epsilon<|x|.
$$
Sending $\epsilon$ to $0$ leads to
$\hat u(x)\ge \delta$, i.e.
$u\ge \delta$ on $\Omega\setminus \{0\}$.
This violates (\ref{ABC6}).
Proposition \ref{prop1} for $k=0$ is established.

\medskip

Now

\noindent\underline{The proof for $k\ge 1$}.\
Let
$$
g(s):=-\ln s + \ln(-\ln s), \qquad 0<s<1,
$$
$$
\hat g(x):= g(d(x)/\bar r), \qquad x\in \Omega\setminus E,
$$
where $d(x):=dist(x, E)$, and $\bar r>9\ diam(\Omega)$
 is some large positive number whose value
will be determined below.
Clearly
$$
g'(s)=-s^{-1}+s^{-1}(\ln s)^{-1},
\qquad
g''(s)=s^{-2}-
s^{-2} (\ln s)^{-1} - s^{-2} (\ln s)^{-2},
$$
and, in view of Lemma \ref{lemapp1}, for some matrix valued function
 $O(x)\in O(n)$,
\begin{eqnarray*}
(\bar r^2)
O(x)^t (\nabla  ^2 \hat g(x)) O(x)
&=& diag\left(0, \cdots, 0, g''(\tilde d),
\tilde d^{-1} g'(\tilde d), \cdots,
\tilde d^{-1} g'(\tilde d)\right)\\
&&+O(1)(\bar r|g'(\tilde d)|+dg''(\tilde d)|),
\end{eqnarray*}
where $\tilde d=d(x)/\bar r$, $|O(1)|\le C$
for some constant $C$ depending only on $E$ --- it is in particular
independent of $\bar r$.
Note that there are $k$ zeros and $n-k-1$ of
$\tilde d^{-1} g'(\tilde d)$ in the above diagonal matrix.
Thus, for $x$ close to $E$,
$$
(\bar r^2)\lambda_i(-\nabla^2 \hat g(x))
=
\left\{
\begin{array}{rl}
-\tilde d^{-2}+ \tilde d^{-2} (\ln \tilde d)^{-1} + \tilde  d^{-2} (\ln
\tilde  d)^{-2}+O(1)\bar r \tilde d^{-1},&
i=1,\\
O(1)\bar r\tilde d^{-1},&
2\le i\le k+1,\\
\tilde d^{-2}- \tilde d^{-2} (\ln \tilde  d)^{-1} + O(1)\bar r\tilde d^{-1},&
k+2\le i\le n,
\end{array}
\right.
$$
and therefore, for some small number $\tilde r\in (0, \frac 19)$ which
depends only on $E$,
\begin{equation}
\sum_{i=1}^{k+2}
\lambda_i( -\nabla^2 \hat g(x))=
 d(x)^{-2}[(\ln \tilde d)^{-2} +O(1)diam(\Omega)]>0,
\qquad \forall\ \tilde d<\tilde r.
\label{AAB2}
\end{equation}
To guarantee $\tilde d<\tilde r$ for all $x$ in
$\Omega\setminus E$, we only need to choose from
the beginning $\bar r> diam (\Omega)/ \tilde r$.

Let
$$
\hat u:= u- \inf_{\Omega\setminus E}u.
$$
Then
\begin{equation}
\inf_{\Omega\setminus E}\hat u=0.
\label{ABC6new}
\end{equation}
We need to prove that
$$
 \inf_{\partial \Omega\setminus E}\hat u=0.
$$
We prove it by contradiction.  Suppose the contrary,  then
\begin{equation}
 \inf_{\partial \Omega\setminus E}\hat u=\delta>0.
\label{ABC4new}
\end{equation}
For
$\bar r= (\tilde r)^{-1} diam (\Omega)$, let
$$
\hat h_\epsilon(x):=
1-\frac {  g\left(\frac{d(x)}{ \bar r}\right)  }
        {   g\left(\frac{\epsilon}{ \bar r}\right)  }.
$$
Then, for $E_\epsilon:= \{x\in \Bbb R^n\ |\ d(x)<\epsilon\}$,
\begin{equation}
\hat h_\epsilon=0\ \ \mbox{on}\ \partial E_\epsilon,\qquad
\hat h_\epsilon<1\ \ \mbox{on}\ \Omega\setminus E_\epsilon.
\label{AABB1}
\end{equation}
By (\ref{ABC1}) and (\ref{AAB2}),
$$
\sum_{i=1}^{k+2}
\lambda_i(\nabla ^2 \hat u)\le 0< \sum_{i=1}^{k+2}
\lambda_i(\nabla ^2 \hat h_\epsilon),\qquad\mbox{in}\
\Omega
\setminus E_\epsilon.
$$
It follows that
$$
\lambda_1(\nabla^2(\hat u-\delta \hat h_\epsilon))<0,\qquad
\mbox{in}\ \Omega
\setminus E_\epsilon.
$$
Therefore, in view of (\ref{ABC6new}), (\ref{ABC4new}) and
(\ref{AABB1}),
$$
\hat u-\delta \hat h_\epsilon\ge
\inf_{ \partial(\Omega
\setminus E_\epsilon) }
(\hat u-\delta \hat h_\epsilon)
\ge 0.
$$
Now, for any $\Omega\setminus E$, we have, from the above, that
$$
\hat u(x)\ge \delta \hat h_\epsilon(x),\qquad \forall\ 0<\epsilon<d(x).
$$
Sending $\epsilon$ to $0$ leads to $\hat u(x)\ge \delta$, i.e.
$u\ge \delta$ on $\Omega\setminus E$.
This violates (\ref{ABC6new}).
Proposition \ref{prop1} is established.

\section{Appendix D. Some properties of viscosity
\newline
$k-$superaffine functions}

We give an extension of Proposition \ref{prop1}.
As always, we
 order the eigenvalues of
a  $C^2$ function $\varphi$  as
$$
\lambda_1(\nabla^2 \varphi)\le \lambda_2(\nabla^2 \varphi)\le
\cdots\le \lambda_n(\nabla^2 \varphi).
$$

Recall that we use ${\cal A}_k(\Omega)$ to denote the set of
viscosity
 $k-$superaffine functions on a bounded open set $\Omega$ of $\Bbb R^n$, and
use $LSC(\Omega)$
to denote the set of lower semi-continuous functions  in $\Omega$.

\begin{thm}
For $n\ge 2$, $0\le k\le n-2$,  let $E$ be a smooth closed
$k-$dimensional manifold in $\Bbb R^n$ and
$\Omega\subset \Bbb R^n$ be a bounded open set.
Assume that
$
u\in
 {\cal A}_{k+2}(\Omega\setminus E).
$
Then, after setting, for $y\in E$,
$
u(y):=\liminf_{x\to y}u(x),
$
$
u\in
 {\cal A}_1(\Omega).
$
\label{Theorem 1.1}
\end{thm}

Theorem \ref{Theorem 1.1}
is somewhat sharp, as explained below.
\begin{definition}
For $n\ge 1$, $0\le k\le n-2$,  $0\le a\le 1$,   let
 $\Omega\subset \Bbb  R^n$ be an open set, and let $u \in
LSC(\Omega)$ satisfy (\ref{DD0}).
We say that $u\in {\cal A}_{k+1,a}(\Omega)$ if:
$$
\varphi\in C^2(\Omega),\
\ \ \varphi \le u,\quad\mbox{in}\ \Omega,\qquad \mbox{and}\ \varphi(\bar
x)=u(\bar x),\ \mbox{for some}\
 \bar x\in \Omega
$$
implies
$$
\lambda_1(\nabla^2 \varphi)+
\cdots+  \lambda_{k+1}(\nabla^2 \varphi)
+ a \lambda_{k+2}(\nabla^2 \varphi)\le
0,\quad \mbox{at}\ \bar x.
$$
\label{def3.1} \end{definition}

 Clearly $${\cal
A}_{k+1,0}(\Omega)={\cal A}_{k+1}(\Omega), \quad  {\cal
A}_{k+1,1}(\Omega)={\cal A}_{k+2}(\Omega),
$$
and
$$
{\cal A}_{k+2}(\Omega)\subset
{\cal A}_{k+1, b}(\Omega)\subset {\cal A}_{k+1, a}(\Omega)\quad \forall\
0\le a\le b\le 1.
$$
For $u\in C^2(\Omega)$, $u\in {\cal A}_{k+1, a}(\Omega)$ means
$$
\lambda_1(\nabla^2 u)+
\cdots+ \lambda_{k+1}(\nabla^2 u)
+ a \lambda_{k+2}(\nabla^2 u)\le 0\ \ \mbox{in}\
\Omega.
$$

\bigskip

\bigskip

The following example shows that in Theorem \ref{Theorem 1.1}
 assumption $u\in
{\cal A}_{k+2}(\Omega\setminus E)$ cannot be replaced by $u\in {\cal
A}_{k+1,a}(\Omega\setminus\{0\})$ for any $0<a<1$.

\begin{exam}

(i)\ Let $n\ge 2$, $k=0$ and $E=\{0\}$.
For any $0<a<1$, let $u(x)=|x|^{1-a}$. Then
$u'(r)=(1-a)r^{-a}>0$, $u''(r)=-a(1-a)r^{-1-a}<0$.  $u'(r)/r=
(1-a)r^{-1-a}$. So, by Lemma \ref{lemapp1}, the
eigenvalues of $\nabla ^2 u$ are
$$
-a(1-a)r^{-1-a}< (1-a)r^{-1-a}=\cdots=(1-a)r^{-1-a}.
$$
So
$$
\lambda_1(\nabla^2 u)+a\lambda_2(\nabla^2u)=0\qquad\mbox{in}\
B_1\setminus\{0\}.
$$
 But $u$ does not belong
to ${\cal A}_1(B)$ since $u(0)=0<u(1)=1$.

(ii)\  For $n\ge 3$, $1\le k\le n-2$,
 we write  $x=(x', x'')\in \Bbb R^n$ with
$x'\in \Bbb R^{n-k}$ and
$x''\in \Bbb R^k$.
Let
$\Omega=B_1$ be the unit ball  in $\Bbb R^n$, and
let $E$ be a $k-$dimensional smooth closed manifold
satisfying $E\cap B_2=\{x=(x', x'')\
|\ x'=0, |x''|<2\}$.
For any $0<a<1$, let $u(x)=|x''|^{1-a-\epsilon}
+\epsilon ^2|x'|^2$ where $\epsilon>0$ is some
small number to be specified below. Then
the  eigenvalues of $\nabla ^2 u$ in $B_1\setminus E$
are
$$
 \lambda_i(\nabla^2 u)=
\left\{
\begin{array}{ll}
 -(a+\epsilon)(1-a-\epsilon)|x''|^{-1-a-\epsilon}&i=1,\\
2\epsilon^2&  2\le i\le k+1,\\
 (1-a-\epsilon)|x''|^{-1-a-\epsilon}& k+2\le i\le n.
\end{array}
\right.
$$
So $u\in C^2(B_1\setminus E)$ and,
for small $\epsilon>0$,
\begin{eqnarray*}&&
\lambda_1(\nabla^2u)+\cdots+\lambda_{k+1}(\nabla^2u)+a\lambda_{k+2}(\nabla^2u)
\\
&=&-\epsilon(1-a-\epsilon)|x''|^{-1-a-\epsilon} +2k\epsilon^2
\le -\epsilon(1-a-\epsilon) +2k\epsilon^2
 <0\qquad \mbox{in}\
B_1\setminus E.
\end{eqnarray*}
But $\displaystyle{
u(0)=0<\epsilon^2=\min_{\partial B_1}u
}$.
\end{exam}

The following two simple
 lemmas give an equivalent definition of
super-affine functions.

\begin{lem} For $n\ge 1$,
let  $\Omega\subset \Bbb R^n$ be an open set,
and let
$u\in  {\cal A}_1(\Omega)$.  Then for every
open subset $D$ satisfying $\overline D\subset \Omega$, and for
every $V\in \Bbb R^n$, $w(x):= u(x)+V\cdot x$ satisfies
\begin{equation}
\min_{ \overline D }w=\min_{\partial D} w.
\label{DD1}
\end{equation}
\label{Lemma 1.1}
\end{lem}

\noindent{\bf Proof.}\  Since $u\in LSC(\Omega)$,
it is clear that both $\inf_{ \overline D }u$ and
$\inf_{\partial D} u$ are achieved.
We prove (\ref{DD1})  by contradiction.  Suppose not, then for some
$D$ and for some $V$,
$$
-\infty<\min_{x\in  \overline D }(u(x)+V\cdot x)<\min_{x\in \partial D}
(u(x)+V\cdot x).
$$
Let
$$
\tilde u(x):= u(x)+V\cdot x.
$$
Then $\tilde u\in  {\cal A}_1(\Omega)$, and
\begin{equation}
\delta:= \min_{\partial D} \tilde u
> \min_{\partial D} \tilde u=, \mbox{say}, 0.
\label{DD2}
\end{equation}
Since $D$ is bounded, there exists some constant
$\bar \epsilon>0$ such that
\begin{equation}
\bar \epsilon |x|^2+\frac \delta 4<
\frac \delta 2,\qquad\forall\ x\in \overline D.
\label{DD3}
\end{equation}
Let
$$
\varphi_a(x):= \bar \epsilon |x|^2+a.
$$
For very negative $a$,
$\varphi_a\le \tilde u$ on $\overline D$.
Let
$$
\bar a:= \sup \{ a\ |\ \varphi_a\le \tilde u\ \mbox{on}\
\overline D\}.
$$
By (\ref{DD2}) and (\ref{DD3}),
$$
\min_{ \overline D}(\tilde u-\varphi_{\bar a})=0,\quad
\bar a\le 0, \quad
\mbox{and}\ \tilde u-\varphi_{\bar a} >\frac \delta 4>0\ \
\mbox{on}\ \partial D.
$$
It follows that for some $\bar x\in D$,
$$
(\tilde u-\varphi_{\bar a})(\bar x)=
\min_{ \overline D}(\tilde u-\varphi_{\bar a})=0.
$$
With the above, we
 can easily find a smooth $\tilde \varphi$ satisfying
$$
(\tilde u- \tilde \varphi)\ge 0\  \mbox{in}\ \Omega,
\ (\tilde u- \tilde \varphi)(x)=0, \
\nabla^2  \tilde \varphi(\bar x)=\nabla ^2\varphi_{\bar a}(\bar x)=
2\bar\epsilon I>0.
$$
Since $u$ is in
${\cal A}_1(\Omega)$, we have
$$
\lambda_1(\nabla^2  \tilde \varphi)(\bar x)\le 0.
$$
This is a contradiction.
Lemma \ref{Lemma 1.1}  is proved.

\bigskip

\begin{lem} For $n\ge 1$, let  $\Omega\subset \Bbb R^n$ be an open set,
and let $u\in LSC(\Omega)$ satisfy
(\ref{DD0}); $u$ is bounded below.
If  for every $\bar x\in \Omega$, there exists
some $0<\bar\epsilon<dist(\bar x, \partial \Omega)$
 such that for every $0<\epsilon<\bar\epsilon$,
 and for
every $V\in \Bbb R^n$, $w(x):= u(x)+V\cdot x$ satisfies
$$
%\begin{equation}
\min_{ \overline B_\epsilon(\bar x) }w=\min_{\partial B_\epsilon(\bar x)} w,
$$
%\label{DD5}
%\end{equation}
then $u\in  {\cal A}_1(\Omega)$.
\label{Lemma 1.2}
\end{lem}

\noindent{\bf Proof.}\
We prove it by contradiction.  Suppose that $u$ is
not in ${\cal A}_1(\Omega)$, then for some $\varphi\in C^2(\Omega)$, some
$\bar x \in \Omega$, we have
$$
\varphi \le u,\ \mbox{in}\ \Omega,
\ \ \varphi(\bar x)=u(\bar x),\ \
\lambda_1(\nabla^2 \varphi(\bar x))>0.
$$
For $x$ close to $\bar x$,
$$
\varphi(x)\ge \varphi(\bar x)+
\nabla \varphi(\bar x)(x-\bar x)+
\frac 12 \lambda_1(\nabla^2 \varphi(\bar x))
|x-\bar x|^2.
$$
Thus, for small $\epsilon>0$, and for $V:=\nabla \varphi(\bar x)$,
$$
u(x)- V\cdot
x
\ge \varphi(x)- V\cdot
x
 > \varphi(\bar x)-V\cdot \bar x,
\qquad\forall\ x\in \partial B_\epsilon(\bar x),
$$
while
$$
u(\bar x)- V\cdot
\bar x=  \varphi(\bar x)-V\cdot \bar x.
$$
This implies
$$
\min_{ x\in \overline B_\epsilon(\bar x) }
(u(x)- V\cdot
x)< \min_{ x\in \partial  B_\epsilon(\bar x) }
(u(x)- V\cdot
x).
$$
A contradiction.

\bigskip

\begin{lem}
For  $1\le k\le n$, if $u\in  {\cal A}_k(\Omega)$ and
$v\in C^2(\Omega)$ satisfies

$$
 \lambda_1(\nabla^2 v)+
\cdots+\lambda_k(\nabla^2 v)\ge 0\ \
\mbox{in}\ \Omega,
$$
then
$$
u-v\in {\cal A}_1(\Omega).
$$
\label{Lemma 1.3}
\end{lem}

\noindent{\bf Proof.}\ Suppose not,
then for some $u\in
{\cal A}_k(\Omega)$ and   some $v\in C^2(\Omega)$ satisfying
$$
\lambda_1(\nabla^2 v)+
\cdots+\lambda_k(\nabla^2 v)\ge 0\ \mbox{in}\ \Omega,
$$
$u-v$ is not in ${\cal A}_1(\Omega)\}$.
 Then for some $\varphi\in C^2(\Omega)$, some
$\bar x \in \Omega$, we have
$$
\varphi \le u-v,\ \mbox{in}\ \Omega,
\ \ \varphi(\bar x)=(u-v)(\bar x),\ \
\lambda_1(\nabla^2 \varphi(\bar x))>0.
$$
Since $u\in {\cal A}_k(\Omega)\}$, it follows from the above that
$$
\lambda_1(\nabla^2 (\varphi+v)(\bar x))
+\cdots+\lambda_k(\nabla^2 (\varphi+v)(\bar x))\le 0.
$$
Since $\nabla^2\varphi(\bar x)$ is positive definite,
we arrive at
$$
\lambda_1(\nabla^2 v(\bar x))
+\cdots+\lambda_k(\nabla^2 v(\bar x))
< 0.
$$
this is a contradiction.

\bigskip

\bigskip

\noindent{  \bf Proof of Theorem \ref{Theorem 1.1}.}\
We prove it by contradiction.  Suppose the contrary, then,
there exists
  $u\in {\cal A}_{k+2}(\Omega\setminus E)$,
but,  after setting
$u(y):=\liminf_{x\to y}u(x),
$
$u$ does not belong to
$
 {\cal A}_1(\Omega).
$
By Lemma \ref{Lemma 1.2},
  there exist $\bar x\in \Omega$,
 $\epsilon_i>0$,  $\epsilon_i\to 0$,
$V_i\in \Bbb R^n$, such that
$w_i:= u(x)+V_i\cdot x$ satisfies
\begin{equation}
\min_{ \overline B_{\epsilon_i}(\bar x) }w_i<
\min_{\partial  B_{\epsilon_i}(\bar x) }w_i.
\label{DD7}
\end{equation}
Since $u\in {\cal A}_{k+2}(\Omega\setminus\{0\})\subset
 {\cal A}_1(\Omega\setminus\{0\})$, we must have
$\bar x=0$.

Let
$$
g(s)=
\left\{
\begin{array}{ll}
-\ln s, & \mbox{if}\ k=0,\\
-\ln s+\ln(-\ln s),& \mbox{if}\ 1\le k\le n-2.
\end{array}
\right.
$$
It was proved in Appendix C, see (\ref{AAB2}),
that for some small positive constant $\beta$
$$
\hat g(x):= g(d(x)/\beta), \qquad x\in \Omega\setminus E
$$
satisfies
$$
\hat g(x)>0, \quad
\sum_{i=1}^{k+2}
\lambda_i( -\nabla^2 \hat g(x))>0,\qquad
\forall \ x\in \Omega\setminus E.
$$

Fixing some $i=\bar i$ for which
(\ref{DD7}) holds, and denote $\bar r=
 \epsilon_{\bar i}$,
and
$$
%\begin{equation}
\hat u:= w_{\bar i}-
 \inf_{B_{\bar r}} w_{\bar i}.
$$
%\label{44}
%\end{equation}
Clearly, $\hat u\in  {\cal A}_{k+2}(B_{\bar r}\setminus E)$,
and
\begin{equation}
\inf_{B_{\bar r}}
\hat u=0,
\label{3}
\end{equation}
$$
%\begin{equation}
\hat u\ge \delta:=
 \min_{ \partial B_{\bar r} } \hat u= \min_{ \partial B_{\bar r} }w_{\bar r}
-\inf_{B_{\bar r}}w_{\bar r}>
 0 \quad\mbox{on}\  \partial B_{\bar r}.
$$
%\label{4}
%\end{equation}

For any $\epsilon>0$ small, let
$$
\hat h_\epsilon(x):=
1-\frac {  g\left(\frac{d(x)}{ \bar r}\right)  }
        {  g\left(\frac{\epsilon}{ \bar r}\right)  }.
$$
Then
(\ref{AABB1}) holds for
$E_\epsilon:= \{x\in \Bbb R^n\ |\ d(x)<\epsilon\}$.

By the property of $\hat g$,
$$
\lambda_1(\nabla^2  \hat h_\epsilon)+
\cdots+\lambda_k(\nabla^2  \hat h_\epsilon)
\ge 0,\qquad\mbox{in}\ B_{\bar r}\setminus \overline E_\epsilon.
$$
Thus, by
 Lemma \ref{Lemma 1.3},
$\hat u-\hat h_\epsilon\in
{\cal A}_1(B_{\bar r}\setminus \overline E_\epsilon)
$
It follows that
$$
\hat u-\delta \hat h_\epsilon\ge
\inf_{ \partial  (B_{\bar r}\setminus E_\epsilon) }
(\hat u-\delta \hat h_\epsilon)
\ge 0.
$$
Sending $\epsilon$ to $0$ leads to $\hat u\ge \delta$  on
$B_{\bar r}\setminus E$.  This violates
(\ref{3}).
Theorem \ref{Theorem 1.1}   is established.

\section{Appendix E. Proof of Lemma \ref{abc}}

\noindent{\bf Proof of Lemma \ref{abc}.}\
Suppose not, then there are $n+1$ vectors $p_0, p_1, \cdots, p_n$
satisfying (\ref{(3)}) with
$p_1-p_0, \cdots, p_n-p_0$  linearly
independent.

After a linear transformation of variables, and adding
a linear function to $u$ --- this preserves
(\ref{superaffine}) --- we may
assume that
$$
p_0=0, p_1=e_1, \cdots, p_n=e_n,
$$
where $e_1=(1, 0, \cdots, 0)$, $e_2=(0, 1, 0, \cdots, 0)$, etc.

So
\begin{equation}
u(x)\ge \circ(|x|)
\ \mbox{and}\ u(x)\ge x_i+\circ(|x|),\
i=1, \cdots, n.
\label{(5)}
\end{equation}

Consider
$$
w=u(x)-\frac 1{2n}\sum_{i=1}^n x_i.
$$

\noindent\underline{Claim}:\
$$
w(x)\ge \frac 1{2n}|x|+\circ(|x|).
$$
This would then contradict (\ref{superaffine}).

\noindent\underline{Proof of Claim}:\ Consider any $x$ near the origin.
Suppose, say
$$
x_1, \cdots, x_k\ \mbox{are}\
positive,
\ \
x_{k+1}, \cdots, x_n\ \mbox{are}\
\le 0.
$$
Then
\begin{eqnarray*}
w(x)&=&
\frac 1n \sum_i (u(x)-
\frac 12 x_i)\\
&\ge & \frac 1{2n}
\sum_{i=1}^k x_i
-\frac 1{2n} \sum_{i=k+1}^n x_i+\circ(|x|),\quad \mbox{by}\ (\ref{(5)})\\
&\ge & \frac 1{2n}|x|+\circ(|x|).
\end{eqnarray*}

\vskip 5pt
\hfill $\Box$
\vskip 5pt

Under the assumption of
Lemma A, see (\ref{(3)}),    either $P=\emptyset$ or $P$ contains only
one element.
However, this would not be true under the assumption
of Lemma \ref{abc}, as shown by the following example.

\begin{exam}
 In $\Bbb R^n$, $n\ge 2$, take
$$
u(x)=\max\{0, x_1\}=\left\{
\begin{array}{ll}
0, & x_1\le 0,\\
x_1, & x_1>0.
\end{array}
\right.
$$
Then $$ P=\{ce_1\ |\ 0\le c\le 1\}, $$ where  $P$ is defined as in
Lemma \ref{abc} and $e_1=(1,0,\cdots, 0)$. Clearly, $u\in {\cal
A}_{n-1}(\Bbb R^n)\subset {\cal A}_1(\Bbb R^n)$.
\end{exam}

 Here is another example.

\begin{exam}
For $1\le k<n$, let
$$
u(x)=\max\{0, x_1, \cdots, x_k\}.
$$
Then $$P=\{c_1e_1+\cdots+c_ke_k\ |\ 0\le c_i\le 1\}.$$ Since $u$ is
a function of $k$ variables only,  $u\in {\cal A}_{n-k}(\Bbb
R^n)\subset {\cal A}_1(\Bbb R^n)$.
\end{exam}

\section{ Appendix F}

Let $\Omega\subset \Bbb R^n$ be a bounded domain, and let
 $u\in L^1_{loc}(\Omega)$ satisfies, for some constant $C$,
\begin{equation}
\Delta u\le C\qquad
\mbox{in}\ \Omega
\label{F1-1-new}
\end{equation}
in the
distribution sense.
Then  $u(y)-\frac C2 |y|^2$ is superharmonic in
$\Omega$,
$\av_{ \partial B_r(x) }[u(y)-\frac C2 |y|^2]dy$ and
$\av_{  B_r(x) }[u(y)-\frac C2 |y|^2]dy$ are 
 nonincreasing functions in $r$ in
$(0, dist(x, \partial \Omega))$.
It follows that
$$
u^*(x):= \lim_{r\to 0^+} \av_{ B_r(x) } u 
=  \lim_{r\to 0^+} \av_{ \partial B_r(x) } u
\ \ \ \in (-\infty, \infty],
\qquad x\in \Omega,
$$
satisfies
$$
u^*(x)=\liminf_{y\to x} u^*(y), \qquad x\in \Omega,
$$
and, by a classical theorem of Lebesgue,
$$
u^*=u,\qquad a.e.\ \mbox{in}\ \Omega.
$$
We call $u^*$ the canonical representative
of $u$ in $\Omega$.

\begin{lem}
Let $E\subset \Omega$ be a closed subset of
zero  capacity.
  Assume that $u\in L^1_{loc}(\Omega\setminus E)$
satisfies, for some constant $C$,
\begin{equation}
\Delta u\le C\qquad\mbox{in}\ \Omega\setminus E\
\mbox{in the distribution sense},
\label{80-1}
\end{equation}
and
$$
\inf_{ \Omega\setminus E }u>-\infty.
$$
Then
$
u\in L^1_{loc}(\Omega),
$
and
$$
\Delta u\le C\qquad\mbox{in}\ \Omega\
\mbox{in the distribution sense}.
$$
\label{lem-4.1-new}
\end{lem}

\noindent{\bf Proof.}\  We may assume
that  $u$ is the
canonical representative in $\Omega\setminus E$.
Subtracting a quadratic polynomial from $u$, we may assume without loss of
generality
that $C=0$ and $u\ge 0$ in
$\Omega\setminus E$.
Since $E$ has zero capacity, there exists
a positive harmonic
function $h$ in $\Omega\setminus E$ satisfying
$$
\lim_{ dist(x, E)\to 0}h(x)=\infty;
$$
see  \cite{Evans}.

For constants $\epsilon>0$ and $M>1$, consider
$$
u_{M, \epsilon}:= \min\{u+\epsilon h, M\}.
$$
Since $u\ge 0$ and $h(x)$ tends to infinity as
$x$ approaches $E$,
 $u_{M, \epsilon}=M$ in
a neighborhood of $E$.
Thus
\begin{equation}
\Delta u_{M, \epsilon}\le 0\ \mbox{in}\ \Omega,
\qquad\mbox{in the distribution sense}, 
\label{81-1}
\end{equation}
and,  
for almost all
$\bar x$ in $\Omega\setminus E$, and
for every $0<r<dist(\bar x, \partial \Omega)$,
$$
u_{M, \epsilon}(\bar x)\ge \av_{B_r(\bar x) }u_{M, \epsilon}.
$$
Sending $M$ to infinity, and then $\epsilon $ to zero, we obtain,
using the Fatou theorem, 
$$
u(\bar x)\ge \av_{ B_r(\bar x) }u,\qquad
a.e. \ \bar x\in \Omega\setminus E, \ 0<r<dist(\bar x, \partial \Omega).
$$
Since $u\in L^1_{loc}(\Omega\setminus E)$, $u$ is finite
a.e. in $\Omega\setminus E$.  Therefore we have
$u\in L^1_{loc}(\Omega)$.
By (\ref{81-1}),
$$
\int_\Omega u_{M, \epsilon}\Delta \varphi\le 0,
\qquad \forall\ \varphi\in
C_c^\infty(\Omega), \ \varphi\ge 0.
$$
Sending $\epsilon$ to zero, then $M$ to infinity, we have,
using the Lebesgue dominated convergence theorem,
$$
\int_\Omega u\Delta \varphi\le 0,
\qquad \forall\ \varphi\in
C_c^\infty(\Omega), \ \varphi\ge 0,
$$
i.e. (\ref{80-1}) for $C=0$.
Lemma \ref{lem-4.1-new} is established.

\vskip 5pt
\hfill $\Box$
\vskip 5pt

\begin{definition}
 $u$ is called quasi-superharmonic in
$\Omega$  if
 $u\in L^1_{loc}(\Omega)$ satisfies (\ref{80-1})
 for some constant $C$.
\end{definition}

In the following, if $u$ is quasi-superharmonic in
$\Omega$, we work with its canonical representative 
in $\Omega$, and still denote it by $u$.

We establish
the following improvement of Lemma A.

\begin{lem}   Let $u$ be quasi-superharmonic in $\Omega$, and
let  $\omega$
be a non-negative non-decreasing continuous function on
$(0, 2d), \ d=diam(\Omega)$.
Assume that
  $u$ satisfies, for some $\bar x, \bar y\in \Omega$ and
$p, q\in \Bbb R^n$,
\begin{equation}
u(y)\ge u(\bar x)+p\cdot (y-\bar x)-|y-\bar x|\omega(|y-\bar x|),\qquad
y\in \Omega,
\label{I1}
\end{equation}
and
\begin{equation}
u(y)\ge u(\bar y)+q\cdot (y-\bar y)-|y-\bar y|\omega(|y-\bar y|),\qquad
y\in \Omega.
\label{I2}
\end{equation}
Then, for some constant $C_1$ depending only
on $n$,
\begin{equation}
|p-q|\le
C_1
\omega(2|\bar x-\bar y|)+C_1|\bar x-\bar y|.
\label{F2-1}
\end{equation}
\label{lem-F2}
\end{lem}

Lemma \ref{lem-F2} follows from the following more general lemma,
 which is used in the proof of
Theorem \ref{thmnew1} (see the proof of Lemma \ref{lem-35}).

\begin{lem}
  Let $u$ be quasi-superharmonic in $\Omega$, and
let   $\omega$ be a 
 non-negative  non-decreasing continuous function  on
$(0, 2d), \ d:=diam(\Omega))$.
Assume that   $u$ satisfies, 
for some $\bar x, \bar y\in \Omega$,
$\bar x \ne \bar y$,  and  $p, q$  in $\Bbb R^n$,
that
$$
%\begin{equation}
u(y)\ge u(\bar x)+p\cdot (y-\bar x)-|y-\bar x|\omega(|y-\bar x|),\qquad
y\in \Omega,
$$
%\label{I3}
%\end{equation}
and
\begin{equation}
u(z)\ge u(\bar x)+p\cdot (\bar y-\bar x)-
|\bar y-\bar x|\omega(|\bar y-\bar x|)+ q\cdot (z-\bar y)
-|z-\bar y|\omega(|z-\bar y|),\ \ 
\forall\ z\in \Omega.
\label{I4}
\end{equation}
Then,  for some positive constants $C_1$ and $C_2$ depending only
on $n$,
$$
%\begin{equation}
|p-q|\le C_1 
\omega(2|\bar x-\bar y|)
+C_2 C|\bar x-\bar y|,
$$
% \label{F3-1}
%\end{equation}
\label{lem-F3}
\end{lem}
where $C$ is the constant in (\ref{F1-1-new}).

\noindent{\bf Proof of Lemma \ref{lem-F3}.}\
Working with $\tilde u(z)= u(z+\bar x)-[u(\bar x)+p\cdot z]$
instead of $u(z)$, we may assume without loss of generality
that $\bar x=0$, $u(0)=0$, $p=0$, $\bar y\ne 0$, $q\ne 0$:
\begin{equation}
u(z)\ge -|z|\omega(|z|), \qquad z\in \Omega,
\label{F4}
\end{equation}
\begin{equation}
u(z)\ge - |\bar y|
\omega(|\bar y|)+ q\cdot (z-\bar y)
-|z-\bar y|\omega(|z-\bar y|),\qquad
\forall\ z\in \Omega.
\label{F5}
\end{equation}
By (\ref{F4}),
\begin{equation}
u(z)\ge - 2|\bar y| \omega(2|\bar y|),
\qquad \forall\ 0<|z|\le 2|\bar y|.
\label{F8}
\end{equation}
For  $|z-\bar y|\le \frac 12 |\bar y|$,
 we deduce from
 (\ref{F5}) that
$$
u(z)\ge
 -
|\bar y|\omega(|\bar y|)
+ q\cdot (z-\bar y)
-\frac 12 |\bar y|
\omega(\frac 12|\bar y|)
\ge
-
2|\bar y|\omega(2|\bar y|)
+ q\cdot (z-\bar y).
$$
It follows that
\begin{equation}
u(z)\ge
 -
2|\bar y|\omega(2|\bar y|) +\frac 14 |q||\bar y|, 
\qquad \forall\ z\in 
(B_{ 2 |\bar y|}(0) \setminus   B_{ |\bar y|/2}(\bar y))\cap U,
\label{F9}
\end{equation}
where 
$$
U:= \{z\in \Bbb R^n\ |\ 
q\cdot (z-\bar y)\ge \frac 12 |q||z-\bar y|\}.
$$
Since $(B_{ 2 |\bar y|}(0) \setminus  
B_{ |\bar y|/2}(\bar y))\cap U$ is a subset of
$B_{2 |\bar y|}(0)$ of measure $\ge |\bar y|^n/C_1$
for some positive dimensional constant $C_1$, and since
$u(x)-\frac C{2n}|x|^2$ is superharmonic by (\ref{F1-1-new}),
 we
deduce from (\ref{F8}) and (\ref{F9}),
 in view of the mean value theorem,
that
$$
0=u(0)=  \av_{  B_{2|\bar y|}(0)} \left(u(x)-\frac C{2n}|x|^2\right)
 \ge  - 2|\bar y| \omega(2|\bar y|)
+ 
\frac 1{C_1} |\bar y| |q|
-C_2 C|\bar y|^2
$$
It follows that
$$
|q|\le C_1\omega(2|\bar y|))
+C_2 C|\bar y|.
$$
Lemma \ref{lem-F3} is established.

\vskip 5pt
\hfill $\Box$
\vskip 5pt

\noindent{\bf Proof of Lemma \ref{lem-F2}.}\
Taking $y=\bar y$ in (\ref{I1}), 
and using this to substitute for $u(\bar y)$
in (\ref{I2}), we have (\ref{I4}).
Thus 
(\ref{F2-1}) in the case $\bar x\ne \bar y$ follows from 
Lemma \ref{lem-F3}.

Now we discuss the case when $\bar x=\bar y$.  Replacing $u(y)$ by
$u(y)-[u(\bar x)+p\cdot (y-\bar x)]$, we may assume without
loss of generality that $\bar x=\bar y=0$, $p=0$,
and $q\ne 0$. So we have
$$
%\begin{equation}
u(y)\ge -|y|\omega(|y|),\qquad y\in \Omega,
$$
%\label{J1}
%\end{equation}
\begin{equation}
u(y)\ge  q\cdot y-|y|\omega(|y|),\qquad y\in \Omega.
\label{J2}
\end{equation}

It is easy to see that for some large
 constant $\alpha>1$ independent of
$\epsilon$,
$$
w_\epsilon(y):= \alpha \omega(\alpha |y|+\alpha \epsilon)
$$
satisfies
$$
|y|\omega(|y|)\le \epsilon |q|
+\epsilon\omega(\epsilon)
+|y-\epsilon q/|q||\omega_\epsilon(|y-\epsilon q/|q||).
$$
It follows from  
(\ref{J2}) and the above that
$$
u(y)\ge -|y_\epsilon|\omega_\epsilon(|y_\epsilon|)
+ q\cdot (y-y_\epsilon)-
|y-y_\epsilon|\omega_\epsilon(|y-y_\epsilon|),\quad y\in \Omega,
$$
where $y_\epsilon = \epsilon q/|q|\ne 0$.  We can apply 
Lemma \ref{lem-F3}, with $\bar y=y_\epsilon$, to obtain,
for some constants $C_1, C_2$ depending only on
$n$, 
$$
|q|\le C_1 \omega_\epsilon(2|y_\epsilon)+C_2C|y_\epsilon|,
$$
where $C$ is the constant 
in  (\ref{F1-1-new}).
Sending $\epsilon$ to zero in the above leads to 
$|q|\le C_1\omega(0)$.  Lemma \ref{lem-F2} is established.

\vskip 5pt
\hfill $\Box$

\end{document}